\documentclass[11pt,abstract=on]{scrartcl}
\usepackage[utf8]{inputenc}

\usepackage{authblk}

\usepackage{dsfont}
\usepackage{amsmath}
\usepackage{amsfonts}
\usepackage{amssymb}
\usepackage{amsthm}
\usepackage{mathtools}
\usepackage{hyperref}
\usepackage{subcaption}
\usepackage{booktabs}
\usepackage{tikz}
\hypersetup{hidelinks}

\usepackage[ruled,vlined]{algorithm2e}
\usepackage{algorithmic}
\SetKwInput{Input}{Input}
\SetKwInput{Output}{Output}

\usepackage[
  a4paper,
  left=2.53cm,
  right=2.53cm,
  top=2.49cm,
  bottom=2.64cm
]{geometry}

\allowdisplaybreaks
\usepackage{xcolor}

\usepackage{natbib}
\numberwithin{equation}{section}
\numberwithin{table}{section}
\numberwithin{figure}{section}

\newcommand{\R}{\mathds{R}}

\newcommand{\N}{\mathds{N}}
\newcommand{\X}{\mathcal{X}}
\newcommand{\Y}{\mathcal{Y}}
\newcommand{\U}{\mathcal{U}}

\newcommand{\proj}{\Lambda^{p-1}}

\newtheorem{theorem}{Theorem}[section]
\newtheorem{lemma}[theorem]{Lemma}
\newtheorem{definition}[theorem]{Definition}

\theoremstyle{remark}
\newtheorem{remark}[theorem]{Remark}

\definecolor{mediumblue}{rgb}{0.0, 0.0, 0.8}
\definecolor{orange-red}{rgb}{1.0, 0.27, 0.0}
\definecolor{darkblue}{rgb}{0.0, 0.0, 0.55}
\definecolor{navyblue}{rgb}{0.0, 0.0, 0.5}
\definecolor{royalblue(web)}{rgb}{0.25, 0.41, 0.88}
\definecolor{royalblue(traditional)}{rgb}{0.0, 0.14, 0.4}
\definecolor{royalpurple}{rgb}{0.47, 0.32, 0.66}
\definecolor{bluegray}{rgb}{0.4, 0.6, 0.8}
\definecolor{forestgreen(web)}{rgb}{0.13, 0.55, 0.13}

\title{Pareto Front Tracing in Higher Dimensions Using Numerical Integration}

\author{Lara L\"ohken\thanks{Department of Mathematics and Informatics, Bergische Universit\"at Wuppertal, Wuppertal, Germany. L.~L\"ohken (corresponding author): \href{mailto:loehken@math.uni-wuppertal.de}{loehken@math.uni-wuppertal.de}, \href{https://orcid.org/0009-0003-9205-6641}{ORCID 0009-0003-9205-6641}; K.~Klamroth: \href{mailto:klamroth@math.uni-wuppertal.de}{klamroth@math.uni-wuppertal.de}, \href{https://orcid.org/0000-0001-9119-0732}{ORCID 0000-0001-9119-0732}.}
,\quad Matthias Bolten\thanks{Faculty of Engineering, Friedrich-Alexander-Universit\"at Erlangen-N\"urnberg, Erlangen, Germany; \href{mailto:matthias.bolten@fau.de}{matthias.bolten@fau.de}; \href{https://orcid.org/0000-0002-8682-7652}{ORCID 0000-0002-8682-7652}.}
,\\ Hanno Gottschalk\thanks{Institute of Mathematics, Technische Universit\"at Berlin, Berlin, Germany; \href{mailto:gottschalk@math.tu-berlin.de}{gottschalk@math.tu-berlin.de}; \href{https://orcid.org/0000-0003-2167-2028}{ORCID 0000-0003-2167-2028}.}
,\quad Kathrin Klamroth\textsuperscript{*}}

\date{}

\begin{document}

\maketitle

\begin{abstract}
We propose a numerical integration approach for tracing the Pareto front of unconstrained continuous multi-objective optimization problems. By stating optimality conditions for weighted sum scalarizations, we derive parametric equations for Pareto critical points, parameterized by the weights. Sensitivity analysis with respect to the weights yields a system of partial differential equations (PDEs). We prove local existence and uniqueness of solutions and derive conditions under which a solution extends to the entire weight space and yields Pareto optimal solutions. A tracing strategy for the resulting PDE system is developed that uses numerical  integration along individual weight coordinates using Runge--Kutta methods. We provide an error analysis considering inexact initial solutions, perturbations of the right-hand sides, and numerical integration errors. We investigate the accuracy and computational efficiency of the tracing strategy on convex quadratic problems and WFG examples, where the latter demonstrate the applicability to problems with Pareto fronts consisting of convex and nonconvex parts. Comparisons with repeated weighted sum scalarizations and a predictor--corrector method illustrate the computational benefits of the proposed strategy. \\


\end{abstract}

\section{Introduction} \label{sec:Intro}

Many real-world optimization problems involve balancing several, often conflicting, objectives. Besides classical economic criteria such as cost or profit, environmental, technical and social aspects 
have to be taken into account. We consider an unconstrained multi-objective optimization problem \eqref{eq:MOP}
\begin{equation}\tag{MOP}\label{eq:MOP}
\begin{array}{ll}
    \min & f(x)=(f_1(x),\dots,f_p(x))\\
    \text{s.t.} & x\in \R^n \; , 
\end{array}\end{equation}
where $f$ is second-order continuously differentiable.

When optimizing multiple objectives simultaneously, we generally expect to find different optimal solutions with respect to the individual objectives. In this situation, we adopt the concept of Pareto optimality. A solution is called Pareto optimal if no objective can be improved without deteriorating at least one other objective. The image of the Pareto optimal set in the objective space is referred to as the Pareto front. See \citep{ehrg:mult:2005,miet:nonl:1998} for an introduction into the field of multi-objective optimization.

A variety of methods exist to determine individual Pareto optimal solutions. Gradient-based methods iteratively converge to locally Pareto optimal solutions, and scalarization techniques transform the multi-objective problem into auxiliary single-objective optimization problems. While these approaches are well suited for computing individual Pareto optimal solutions, accurately reproducing the entire Pareto front typically requires to repeatedly solve optimization problems for different initializations or scalarization parameters. This process often becomes computationally demanding, especially for high-dimensional or computationally expensive optimization problems. Repeatedly solving scalarizations with varying parameters hence suggests to investigate the dependence of Pareto optimal solutions on these parameters. 

Following upon this idea, continuation methods successively determine nearby solutions as the parameters vary, starting from a given initial solution, see, e.\,g., \citep{allg:intr:2003,hill:nonl:2001}. Implementations of continuation methods are often based on predictor--corrector methods. In the predictor step, a tangent direction at the current solution is computed and, subsequently, it is corrected by projections onto the Pareto set in a corrector step~\citep{mart:onco:2016,mart:predcorr:2018}. A more detailed review of related literature is provided in Section~\ref{subsec:litrev}.

In this work, we propose a different numerical approach for local and global Pareto front tracing. Implicit differentiation of first-order optimality conditions of weighted sum scalarizations of \eqref{eq:MOP} with respect to the weights yields a system of differential equations. Given an initial Pareto optimal solution, rather than successively generating nearby solutions with predictor--corrector methods, we solve the resulting differential equation system directly by numerical integration and thereby trace the Pareto front. For the bi-objective case, a direct numerical integration approach has already been proposed by~\cite{bolt:trac:2021}. Extending this idea to an arbitrary number of objectives, however, fundamentally changes the mathematical framework. While in the bi-objective case implicit differentiation with respect to a single parameter yields an ordinary differential equation, the general multi-objective case requires differentiation with respect to several independent parameters, resulting in a system of partial differential equations. This, in turn, requires new theoretical and numerical developments.

The main contributions of this paper are summarized as follows. We derive a framework for tracing the Pareto front of an unconstrained multi-objective optimization problem with an arbitrary number of objectives by numerical integration. We state assumptions under which the resulting system of differential equations has a unique local solution and derive conditions under which the complete weight space, and consequently the Pareto front, can be traced. We analyze perturbations resulting from inexact initial Pareto optimal solutions, develop an efficient numerical integration scheme for numerically solving the resulting system of differential equations, and derive corresponding error bounds. We demonstrate the accuracy and computational efficiency of the proposed framework on several benchmark problems and compare it with iteratively repeated weighted sum scalarization and another continuation method, namely the Pareto Tracer by \cite{mart:predcorr:2018}. To support reproducible research, the implementation is made publicly available in the associated repository \citep{loehken:repository:2026}.

The remainder of this paper is organized as follows. In Section~\ref{sec:MathBackground}, we introduce the mathematical background and review related work. In Section~\ref{sec:ParetoTracing}, we develop the theoretical framework, followed by the numerical integration scheme and its error analysis in Section~\ref{sec:NumIntegration}. Numerical experiments demonstrating the accuracy and computational efficiency of the approach are presented in Section~\ref{sec:NumExperiments}. We conclude in Section~\ref{sec:Conclusion}.

\section{Mathematical Background and Related Work} \label{sec:MathBackground}

\subsection{Multi-objective Optimization}

We begin by formalizing basic concepts and notation in multi-objective optimization.
\begin{definition}\label{def:ParetoOpt}
Let $\hat{x}, \bar{x} \in \R^n$. We say that $\hat{x}$ \textit{dominates} $\bar{x}$ if $f_i(\hat{x}) \leq f_i(\bar{x})$ for $i = 1, \dots, p$ and $f_i(\hat{x}) < f_i(\bar{x})$ for at least one $i \in \{1, \dots, p\}$.
The solution $\hat{x}$ is called \textit{Pareto optimal} or \textit{efficient} if there is no $x \in \R^n$ that dominates $\hat{x}$. The corresponding image vector $f(\hat{x})$ is called \textit{non-dominated}.
\end{definition} 

We denote the set of efficient solutions in the decision space by $\X_E$ and the set of non-dominated points in the objective space by $\Y_N = f(\X_E)$. We also refer to $\Y_N$ as the Pareto front ($\mathcal{PF}$). Throughout this paper we assume that $\X_E\neq\emptyset$, and that the objective functions of \eqref{eq:MOP} are conflicting in the sense that there exists no solution $x\in\R^n$ that minimizes all objective functions simultaneously.

The notions of local Pareto optimality and local non-dominance are defined analogously, by replacing the feasible set $\R^n$ with an open neighborhood $\U(\hat{x}) \subseteq \R^n$. A necessary condition for a solution $\hat{x}$ to be locally Pareto optimal is that for all directions $d \in \R^{n}$ there exists at least one index $j \in \{1, \dots, p\}$ such that $\nabla f_j(\hat{x})^\top d \geq 0$, i.\,e., there exists no common descent direction for all objectives $f_i$, $i = 1, \dots, p$. In this case, we call $\hat{x}$ \textit{Pareto critical}. Note that while criticality is a necessary optimality condition, it is in general no sufficient condition for $\hat{x}$ to be Pareto optimal.

Depending on the underlying problem, the Pareto front may have different properties and shapes. We say that $\mathcal{PF}$ is convex if the outcome set $\Y = \{ f(x) \colon x \in \R^n \}$ is $\R_+$-convex, i.\,e.\ if the set $\{ f(x) +  s \colon x \in \R^n, s \in \R^p, s_i \geq 0, i =1, \dots, p\}$ is convex. Otherwise, we say that $\mathcal{PF}$ is non-convex. Moreover, we distinguish cases where the Pareto front consists of convex and non-convex parts or is entirely non-convex. Apart from the convexity, we are also concerned with whether the Pareto front is connected or disconnected. In the literature, we understand a disconnected Pareto front as also being discontinuous. We refer to~\citep{ehrg:mult:2005,hub:revtoolkit:2006} for further reading. Indeed, given the above stated unconstrained~\eqref{eq:MOP} along with the condition that the objectives $f_i, \, i=1, \dots, p$ are convex, $\mathcal{PF}$ is connected and convex~\citep{ehrg:mult:2005}. Note that this is a sufficient, but not necessary condition.

\subsection{Weighted Sum Scalarization} \label{subsec:WS}

Weighted sum scalarizations are a common tool to identify Pareto optimal solutions~\citep{ehrg:mult:2005}. Thereby, the multi-objective problem is transformed into parameter-dependent single-objective optimization problems
\begin{equation}\tag{\ensuremath{\mathrm{WS}_{\lambda}}}\label{eq:WS}
\begin{array}{ll} 
    \min & f_{\lambda}(x)  = \langle \lambda, f(x) \rangle =  \displaystyle  \sum_{i=1}^p \lambda_i f_i(x) \\
    \text{s.t.} & x \in \R^n \;, 
\end{array}
\end{equation}
where $\lambda \in \Lambda = \{ \lambda \in \R^p : \lambda_i > 0 \text{ for } i \in \{1, \dots, p \}, \, \sum_{i=1}^p \lambda_i = 1 \} $. The set of meaningful weights $\Lambda$ is also known as the (normalized) weight space. The vector $\lambda$ can be interpreted as a preference parameter that assigns a weight, representing, e.\,g., relative importance, to each of the objectives. Note that the dimension of $\Lambda$ is $\dim(\Lambda)=p-1$ due to the normalization of the weights. Throughout this paper, we set wlog $\lambda_p = 1 - \sum_{i = 1}^{p-1} \lambda_{i}$ and define the projection of the normalized weight space onto its first $p-1$ components as $\Lambda^{p-1} \coloneqq  \{ \lambda \in \R^{p-1} : \lambda_i > 0 \text{ for } i \in \{1, \dots, p-1 \}, \, \sum_{i=1}^{p-1} \lambda_i < 1 \}$. For $\lambda\in\Lambda^{p-1}$, the weighted objective of \eqref{eq:WS} can equivalently be written as $f_\lambda(x) = \sum_{i=1}^{p-1} \lambda_i f_i(x) + \left( 1-\sum_{i=1}^{p-1} \lambda_i \right)f_p(x)$.

\subsubsection{Correctness and Completeness}
In multi-objective optimization we are generally interested in the \textit{correctness} and \textit{completeness} of the chosen scalarization method. While correctness is concerned with the question whether a generated solution is actually Pareto optimal, completeness addresses the question whether the whole non-dominated set $\mathcal{Y}_N$ can be determined by the considered scalarization method by varying the scalarization parameters. Indeed, the weighted sum scalarization is correct for all $\lambda \in \Lambda$, i.\,e.\ for all $\lambda \in \Lambda$, a solution $\hat{x}$ that is optimal for \eqref{eq:WS} is efficient~\citep{ehrg:mult:2005}. However, only convex parts of the Pareto front can be determined via the weighted sum method in general, see, e.\,g.,~\citep{ehrg:mult:2005,das:aclo:1997}. Hence, the weighted sum method guarantees completeness only if $\mathcal{PF}$ of \eqref{eq:MOP} is convex. Note also that even in this case, completeness can only be guaranteed when allowing weights $\lambda\in\text{cl}(\Lambda)$. Despite this shortcoming, the weighted sum method offers two major advantages that we exploit in this work. First, its formulation naturally extends from the bi-objective to the general multi-objective case by simply introducing additional weight parameters. Second, unlike several other scalarizations, it does not introduce additional constraints to the resulting optimization problem. 

\subsubsection{Optimality Conditions for Weighted Sum Scalarizations}
We now state first- and second-order optimality conditions for \eqref{eq:WS}. A solution $\hat{x} \in \R^{n}$ is called (locally) optimal for $f_{\lambda}$ with respect to the parameter $\lambda \in \Lambda$ if $f_{\lambda}(\hat{x}) \leq f_{\lambda}(x)$ for all $x \in \R^{n}$ (for all $x\in \mathcal{U}(\hat{x}))$. We call $\hat{x}$ critical for $f_{\lambda}$ if $\nabla_{x} f_{\lambda}(\hat{x}) = \sum_{i=1}^p \lambda_i \nabla_{x}f_i (\hat{x}) = 0$, which we refer to as the first-order optimality conditions for \eqref{eq:WS}. Let $\nabla^{2}_{x} f_{i}(x)$ be the Hessian matrix of $f_i$ at $x$, $i \in \{1, \dots, p\}$ and likewise $\nabla^{2}_{x} f_{\lambda}(x)$ the Hessian matrix of $f_{\lambda}$. We say that $\hat{x}$ satisfies the strict second-order optimality conditions for $f_{\lambda}$ if $\hat{x}$ is $f_{\lambda}$-critical and $\nabla^{2}_{x} f_{\lambda}(\hat{x})$ is positive definite. If a solution $\hat{x} \in \R^{n}$ satisfies strict second-order $f_{\lambda}$-optimality conditions, then $\hat{x}$ is a strict local minimizer of $f_{\lambda}$ and hence locally Pareto optimal~\citep{hill:nonl:2001,miet:nonl:1998}.

\subsection{Continuation Methods and Related Work}\label{subsec:litrev}
As weighted sum scalarizations define a family of parameter-dependent optimization problems, there is a close relation to the field of parametric optimization; see, e.\,g.,~\citep{neme:survparopt:2025}. The relation between the weight space $\Lambda$ and the resulting distribution of non-dominated points in the objective space is generally nontrivial and problem dependent~\citep{das:aclo:1997}. Therefore, various parameter-selection and scalarization methods have been proposed to improve the distribution of the generated non-dominated points~\citep{Eich09,pere:eqparfr:2013,przy:dicho:2019}. 

Instead of repeatedly solving parameterized optimization problems, which can become computationally demanding, continuation methods exploit the smooth dependence of optimal solutions on the underlying parameters. We refer to~\citep{gudd:para:1987,allg:intr:2003} for an introduction to continuation methods. Many continuation methods use predictor--corrector schemes. They have been applied, e.g., to 
aerodynamic~\citep{schm:pare:2008} and multi-objective shape optimization~\citep{cesa:tracingelectric:2024,cesa:shapenewton:2025}.

\cite{rako:mocontrolhomotopy:1993} extended predictor--corrector continuation strategies for parameter-dependent nonlinear optimization problems~\citep{rao:1989:contstrat,rako:activesetparopt:1991} to constrained bi-objective optimization using active-set techniques. \cite{hill:nonl:2001} proposed a homotopy framework based on the KKT first-order optimality conditions for smooth multi-objective optimization problems. Under suitable regularity assumptions, the resulting solution set locally forms a manifold of Pareto-critical solutions.

Since continuation methods are local by nature, subsequent research followed two complementary directions. On the one hand, continuation was embedded into global optimization frameworks, including set oriented methods~\citep{schu:seto:2013}, subdivision techniques~\citep{dell:cove:2005,jahn:mult:2006} or evolutionary algorithms~\citep{coel:evoalg:2007,deb:mopevo:2001}. 

\cite{schu:onco:2005} embedded a predictor--corrector method based on the continuation framework of~\cite{hill:nonl:2001} into a global subdivision algorithm. 
In the context of evolutionary algorithms, continuation has been combined with particle swarm optimization~\citep{schu:evostratcont:2008}, genetic algorithms~\citep{Hara:locsearch:2006,hara:unisamp:2007}, and NSGA-II~\citep{cuat:hybevo:2020}. Derivative-free continuation strategies based on approximated trade-off directions have furthermore been proposed~\citep{schu:evocontmeth:2009}.

On the other hand, local continuation methods were further developed to improve their applicability and computational efficiency. 
\cite{pere:eqparfr:2013} proposed an arc-length homotopy approach based on the first-order optimality conditions of weighted sum scalarizations of bi-objective constrained problems.

\cite{mart:onco:2016} suggest a method based on the KKT first-order optimality conditions of the multi-objective optimization problem and combines interval-based continuation with active-set management. Building upon this approach,~\cite{mart:predcorr:2018} introduced the Pareto Tracer, a predictor--corrector method that performs continuation along the set of local Pareto critical solutions and is able to handle equality and box constraints. 
\cite{belt:ParTraceIneqConstr:2020} extended the Pareto Tracer to general inequality constraints by adapting active-set techniques. 

Further developments focus on problem classes that pose additional challenges for continuation methods. \cite{ring:hand:2012} consider high-dimensional decision spaces, the Pareto Explorer addresses many-objective optimization~\citep{schu:manyobj:2020}, and \cite{schu:degmop:2025} extend the Pareto Tracer to degenerate multi-objective optimization problems, where the Pareto set may have dimension smaller than $p-1$. \cite{cesa:shapenewton:2025} introduces higher-order predictors by coupling homotopy techniques with shape-Newton methods, thereby further advancing predictor--corrector continuation methods in multi-objective shape optimization.

A different line of research computes continuation paths directly by numerical integration. \cite{pots:trac:2011} derive an ordinary differential equation from normal-boundary intersection (NBI) scalarizations for constrained bi-objective optimization problems. 
\cite{bolt:trac:2021} derive an ordinary differential equation from first-order optimality conditions of weighted sum scalarizations for smooth unconstrained bi-objective optimization problems. Starting from an initial Pareto optimal solution,  numerical integration is used to trace the Pareto front. 
In this work, we extend this approach to general multi-objective problems. To the best of our knowledge, this is the first framework of this kind for arbitrary dimensions. In contrast to the bi-objective case, where continuation is described by a single ordinary differential equation, considering multiple weight parameters yields a system of partial differential equations. Consequently, the extension requires fundamentally new theoretical and numerical developments.

\section{Pareto Front Tracing for an Arbitrary Number of Objectives} \label{sec:ParetoTracing}

By implicitly differentiating first-order optimality conditions of weighted sum scalarizations with respect to the weight parameters, we derive a system of partial differential equations describing the dependence of Pareto critical solutions on the weight parameters. We investigate the local existence and uniqueness of a function that solves this system and discuss the Pareto optimality of the resulting solutions. Finally, we derive conditions under which the local solution can be extended to the complete weight space.

To derive the PDE system with respect to independent weight parameters, recall that, for normalized weight vectors $\lambda\in\Lambda$ introduced in Section~\ref{subsec:WS}, we have $\lambda_p = 1-\sum_{i=1}^{p-1}\lambda_i$. 
In the following, we therefore work with the corresponding projected weight space $\proj$. We first consider multi-objective optimization problems with a convex Pareto front. This allows us to investigate the theoretical properties of the proposed Pareto tracing framework independently of the limitations of the weighted sum method itself. 
The applicability of the proposed framework to problems whose Pareto fronts contain both convex and non-convex parts is investigated in the numerical experiments.

\subsection{System of PDEs for Pareto Critical Points}

\paragraph{Derivation of the PDE System}
Solving \eqref{eq:WS} for $\lambda \in \proj$, any optimal solution must satisfy the first-order optimality conditions, i.\,e.\ $\nabla_x f_\lambda(x) = 0$. We start by assuming that, for each weight vector $\lambda$ in a subset $\mathcal{U}\subseteq\proj$, a corresponding $f_\lambda$-critical solution $x(\lambda)$ is given. We regard $x$ as a function of the $p-1$ independent weight parameters and write $x(\lambda)$ instead of $x(\lambda_1,\dots,\lambda_{p-1})$ for notational convenience. The main idea is to investigate the sensitivity of $x(\lambda)$ with respect to the weight parameters while still satisfying the first-order optimality conditions. We therefore implicitly differentiate the system $\nabla_x f_\lambda(x(\lambda))=0$ with respect to $\lambda_i$, $i=1,\dots,p-1$, using the chain rule, and obtain
\begin{equation}\label{eq:first_order_opt_cond1}
\begin{array}{l} 
    \nabla_x^2f_{\lambda}(x(\lambda)) \, \partial_{\lambda_i}x(\lambda) =  \nabla_x f_p(x(\lambda)) - \nabla_x f_i(x(\lambda)) \quad \forall \; i = 1, \dots, p\!-\! 1.
\end{array} \end{equation}
In contrast to the bi-objective case, where we get an implicit ordinary differential equation (ODE)~\citep{bolt:trac:2021}, the general multi-objective case yields a system of $p-1$ first order partial differential equations (PDEs). The reason for this is the dimension of $\proj$. If the Hessian $\nabla_x^2f_{\lambda}(x(\lambda)) $ is invertible, we can rearrange the system of PDEs
\begin{equation}\label{eq:first_order_opt_cond2}
\begin{array}{l} 
    \partial_{\lambda_i}x(\lambda) =  - (\nabla_x^2f_{\lambda}(x(\lambda)) )^{-1} \, (\nabla_x f_i(x(\lambda)) - \nabla_x f_p(x(\lambda)) ) \quad \forall \; i = 1, \dots, p\!-\! 1.
\end{array} \end{equation}
Now, instead of assuming from the outset that Pareto optimal solutions $x(\lambda)$ are given, we proceed vice versa and start with the system of PDEs and an initial efficient solution $x^0=x(\lambda^0)$ to determine solutions $x(\lambda)$ that are (locally) optimal for \eqref{eq:WS}. The main question that arises is whether a solution to the system of PDEs exists.

\paragraph{Local Existence}

In the bi-objective case, existence and uniqueness of a solution to~\eqref{eq:first_order_opt_cond2} follow from classical ODE theory, in particular the Picard--Lindelöf theorem~\citep[Theorem~15.2]{agar:introode:2008}. In the general multi-objective setting, however, this approach is no longer directly applicable. Although each of the equations in~\eqref{eq:first_order_opt_cond2} for $i = 1, \dots, p-1$ can, to some extent, be interpreted as an ODE, as we will discuss again in Section~\ref{sec:NumIntegration}, the existence of a solution has to be established for the coupled system as a whole. Therefore, the arguments from the bi-objective case cannot be applied directly. Instead, we employ the implicit function theorem~\citep{krantz:implicit:2013} to derive the local existence and uniqueness of a continuously differentiable function satisfying the first-order optimality conditions and show that this function solves the PDE system. Note that, for unconstrained multi-objective optimization problems, the first-order optimality conditions of weighted sum scalarizations coincide with the corresponding KKT conditions. Hence, the following local existence result is analogous to the local continuation framework proposed by Hillermeier~\citep{hill:nonl:2001}.

\begin{theorem}\label{thm:pareto_tracing_neighborhood} 
Let a Pareto optimal solution $x^0 \in \R^n$ and a weight vector $\lambda^0 \in \proj$ be given, such that $x^0$ satisfies the strict second-order optimality conditions for $f_{\lambda^0}$. Then there exist open neighborhoods $\U(\lambda^0)\subseteq\proj$ and $\mathcal{V}(x^0)\subseteq\R^n$ and a unique continuously differentiable function 
$x \colon \U(\lambda^0)\to\mathcal{V}(x^0)$
that satisfies the first-order optimality conditions for all $\lambda\in\U(\lambda^0)$. Moreover, $x$ is the unique solution to \eqref{eq:first_order_opt_cond2} with the initial condition $x(\lambda^0) = x^0$. 
\end{theorem}

\begin{proof}
Define $ F \colon \R^n\times\proj\to\R^n$, $F(x,\lambda):=\nabla_x f_\lambda(x)$. Since $x^0$ satisfies the strict second-order optimality conditions for $f_{\lambda^0}$, it holds that $F(x^0,\lambda^0)=0$. Moreover, the Hessian $\nabla_x^2 f_{\lambda^0}(x^0) = \nabla_xF(x^0,\lambda^0)$ is positive definite and therefore invertible. Hence, by the implicit function theorem, there exist open neighborhoods $\mathcal{U}(\lambda^0)\subseteq\proj$ of $\lambda^0$ and $\mathcal{V}(x^0)\subseteq\R^n$ of $x^0$ and a unique continuously differentiable function $ x_F \colon \mathcal{U}(\lambda^0)\to\mathcal{V}(x^0)$ with $x_F(\lambda^0)=x^0$ such that, for all
$(x,\lambda)\in
\mathcal{V}(x^0)\times\mathcal{U}(\lambda^0)$,
\[
    F(x,\lambda)=0
    \quad\Longleftrightarrow\quad
    x=x_F(\lambda).
\]
Thus,
\[
    \nabla_x f_\lambda(x_F(\lambda))=0
    \qquad
    \forall\,\lambda\in\mathcal{U}(\lambda^0).
\]
Implicit differentiation with respect to $\lambda_i$, $i=1,\dots,p-1$, yields \eqref{eq:first_order_opt_cond2}. To show uniqueness of the solution to \eqref{eq:first_order_opt_cond2}, let
$\hat{x}_F:\U(\lambda^0)\to\mathcal{V}(x^0)$ be another solution satisfying $\hat{x}_F(\lambda^0)=x^0$. Define $G(\lambda) \coloneqq \nabla_x f_\lambda (\hat{x}_F(\lambda))$. We first show that $G(\lambda)$ is constant on $\U(\lambda^0)$. Differentiating of $G(\lambda)$ yields
$$
\partial_{\lambda_i} G(\lambda) = \nabla_x^2 f_\lambda(\hat{x}_F(\lambda)) \, \partial_{\lambda_i} \hat{x}_F(\lambda) + \nabla_x \, f_i(\hat{x}_F(\lambda)) - \nabla_x \, f_p(\hat{x}_F(\lambda)).
$$
Since $\hat{x}_F(\lambda)$ is a solution to \eqref{eq:first_order_opt_cond2}, it follows directly that $\partial_{\lambda_i} G(\lambda) = 0$ for all $i = 1, \dots, p-1$ and $\lambda \in \U(\lambda^0)$. Thus, $G(\lambda)$ is constant on $\U(\lambda^0)$. Since $\hat{x}_F(\lambda^0) = x^0$ and $\nabla_x f_{\lambda^0}(x^0) = 0 $, we obtain $\nabla_x f_\lambda (\hat{x}_F(\lambda)) = 0$ for all $\lambda \in \U(\lambda^0)$. By the uniqueness property of $x_F$ derived above, it follows that $\hat{x}_F(\lambda)=x_F(\lambda)$ for all $\lambda\in\U(\lambda^0)$. Hence, $x_F$ is the unique solution to~\eqref{eq:first_order_opt_cond2} in $\U(\lambda^0)$ with initial condition $x_F(\lambda^0)=x^0$.
\end{proof}

While Pareto optimality does not follow from first-order optimality conditions alone, the assumptions on the initial solution in Theorem~\ref{thm:pareto_tracing_neighborhood} guarantee local Pareto optimality of $x(\lambda)$ as well. Since $f_i\in C^2$ for $i=1,\dots,p$, the Hessian $\nabla_x^2 f_\lambda(x)$ depends continuously on $(x,\lambda)$. Hence, since $\nabla_x^2 f_{\lambda^0}(x^0)$ is positive definite, $\U(\lambda^0)$ can be chosen sufficiently small such that $\nabla_x^2 f_\lambda(x(\lambda))$ is positive definite for all $\lambda\in\U(\lambda^0)$. Together with the first-order optimality conditions, $x(\lambda)$ satisfies the strict second-order optimality conditions for $f_\lambda$ and is locally Pareto optimal for all $\lambda\in\U(\lambda^0)$.

Indeed, the question arises whether the local result can be extended to the complete projected weight space $\proj$.

\paragraph{Global Pareto Tracing}

For the following result, we assume that $f_\lambda$ is strictly convex, such that every stationary point is its unique global minimizer. A sufficient condition for strict convexity of $f_\lambda$ is that the objective functions $f_i$, $i=1,\dots,p$, are strictly convex. In addition, we assume that $f_\lambda$ attains its minimum for every $\lambda\in\proj$. If this is not the case, the considered weight domain can instead be restricted to a subset of $\proj$ on which the corresponding weighted sum problems have optimal solutions.

\begin{theorem} \label{thm:pareto_tracing}
Assume that, for every $\lambda\in\proj$, the weighted objective $f_\lambda$ is strictly convex and attains a minimum at
\[
x(\lambda) = \arg\min_{x\in\R^n} f_\lambda(x).
\]
Let furthermore $\nabla_x^2 f_\lambda(x(\lambda))$ be invertible for all $\lambda\in\proj$. Let $\lambda^0\in\proj$ be arbitrary and set $x^0\coloneqq x(\lambda^0)$. Then the mapping $x \colon \proj\to\R^n$
is continuously differentiable and is the unique solution to \eqref{eq:first_order_opt_cond2} on $\proj$ with initial condition $x(\lambda^0)=x^0$. Moreover, $x(\lambda)$ is Pareto optimal for \eqref{eq:MOP} for all $\lambda\in\proj$.

\end{theorem}

\begin{proof}
As $f_\lambda$ is strictly convex and has a minimum for every $\lambda\in\proj$, this minimum is unique. Thus, the mapping $ x \colon\proj\to\R^n$, $\lambda\mapsto x(\lambda)= \arg\min_{x\in\R^n} f_\lambda(x)$, is well-defined. It remains to show that $x$ is continuously differentiable.

Let $\lambda^k\in\proj$ be arbitrary and set $x^k:=x(\lambda^k)$. By the definition of $x(\lambda^k)$, it holds that $\nabla_x f_{\lambda^k}(x^k)=0$. Since $f_\lambda$ is strictly convex and $\nabla_x^2 f_{\lambda^k}(x^k)$ is invertible, $x^k$ satisfies second-order optimality strictly. Using Theorem~\ref{thm:pareto_tracing_neighborhood} then yields neighborhoods $\U(\lambda^k)$ of $\lambda^k$ and $\mathcal{V}(x^k)$ of $x^k$ and a unique continuously differentiable function $x^k_F:\U(\lambda^k)\to \mathcal{V}(x^k)$ such that, for all $(x,\lambda)\in \mathcal{V}(x^k)\times \U(\lambda^k)$,
$$
\nabla_x f_\lambda(x) = 0 
\quad\Longleftrightarrow\quad
x=x^k_F(\lambda).
$$
For every $\lambda\in \U(\lambda^k)\cap\proj$, we therefore have $\nabla_x f_\lambda(x^k_F(\lambda))=0$. Since $f_\lambda$ is strictly convex, every stationary point of $f_\lambda$ is its unique global minimizer. Consequently,
$$
x^k_F(\lambda) = \arg\min_{x\in\R^n} f_\lambda(x)=x(\lambda)
$$
for every $\lambda\in \U(\lambda^k)\cap\proj$. Thus, $x$ coincides locally with the continuously differentiable function $x_F^k$. Since $\lambda^k\in\proj$ was arbitrary, it follows that $x\in C^1(\proj,\R^n)$. Moreover, $x$ is the unique mapping for which the first-order optimality conditions hold for every $\lambda\in\proj$, and $x$ solves \eqref{eq:first_order_opt_cond2} on $\proj$. 

Uniqueness of $x$ as a solution to \eqref{eq:first_order_opt_cond2} with $x(\lambda^0)=x^0$ follows as in the proof of Theorem~\ref{thm:pareto_tracing_neighborhood}.

By the correctness of weighted sum scalarizations, $x(\lambda)$ is Pareto optimal for \eqref{eq:MOP} for all $\lambda\in\proj$.
\end{proof}

Theorem~\ref{thm:pareto_tracing} shows that the Pareto optimal solutions can be traced over the complete projected weight space $\proj$, and thus $\Lambda$. Hence, under the stated assumptions, the tracing framework does not introduce any additional limitation compared with solving weighted sum scalarizations individually for all $\lambda \in \proj$.

\subsection{Numerical Stability and Dependence on Initial Conditions}\label{subsec:InitCond}

To prove the existence and uniqueness of a solution $x(\lambda)$ of the PDE system stated in \eqref{eq:first_order_opt_cond2}, we assume that we are given an initial Pareto optimal solution $x^0$ that satisfies second-order optimality strictly. Pareto optimal solutions and thus candidates for an initial value can be computed by solving single-objective weighted sum scalarizations \eqref{eq:WS}. However, since for all $\lambda \in \proj$ the problem \eqref{eq:WS} is a continuous optimization problem, it may be computationally challenging -- or even impossible -- to obtain an exact solution. In such cases, we may have to settle for sub-optimal and/or approximate solutions of \eqref{eq:WS}. As we cannot expect that such solutions satisfy the first-order optimality condition for $f_\lambda$, we introduce relaxed versions of criticality. We call a solution $\widetilde{x} \in \R^n$ $\varepsilon$-Pareto critical, if there exists no search direction $d \in \R^n$ such that $\nabla_x f_i(\widetilde{x})^{\top} d < 0$ for $i = 1, \dots, p$ and $\nabla_x f_j(\widetilde{x})^{\top}d \leq -\varepsilon' \lVert d \rVert$ for $j \in \{1, \dots, p \}$, where $\varepsilon' > 0$. If for \eqref{eq:WS} with $\lambda \in \proj$ holds that $\lVert \nabla_x f_\lambda(\widetilde{x}) \rVert \leq \varepsilon$ for some $\varepsilon > 0$, the point $\widetilde{x}$ is called $\varepsilon$-critical with respect to $f_\lambda$. It directly follows that an $\varepsilon$-critical solution is also $\varepsilon'$-Pareto critical with $\varepsilon' = {\varepsilon} /  {\min\{\lambda_1, \dots, \lambda_p\}}$.

In the following, we examine the effect of an input error caused by a sub-optimal initial value $\widetilde{x}^0$ when (numerically) solving the PDE system \eqref{eq:first_order_opt_cond2}. There exist several optimization algorithms such as gradient descent or Newton's method, which are widely used to approximate (local) minimizers of continuous optimization problems~\citep{noce:numopt:2006}. Starting from an initial value sufficiently close to $x^0$, these methods generate a sequence of iterates $(x_m^0)_{m \in \N}$ that converges to $x^0$ under suitable regularity assumptions. A common stopping criterion of such algorithms is $\lVert \nabla_x f_{\lambda^0}(x_m^0) \rVert \leq \varepsilon$ for some $\varepsilon > 0$. Since $x_m^0 \to x^0$ and $\nabla_x^2 f_{\lambda^0}(x^0)$ is positive definite, we can choose $m^* \in \N$ sufficiently large such that both $\lVert \nabla_x f_{\lambda^0}(x_{m^*}^0) \rVert \leq \varepsilon$ and $\nabla_x^2 f_{\lambda^0}(x_{m^*}^0)$ is positive definite. Define $\widetilde{x}^0 \coloneqq x_{m^*}^0$. By definition, $\widetilde{x}^0$ is $\varepsilon$-critical and therefore also $\varepsilon'$-Pareto critical. 

\begin{theorem}\label{thm:pareto_tracing_init_perturbed_local}
    Let $x^0 \in \R^n$ be a Pareto optimal solution that satisfies second-order optimality strictly for $f_{\lambda^0}$ with respect to $\lambda^0 \in \proj$, let $(x_m^0)_{m \in \N}$ be a sequence of iterates converging to $x^0$, and let $\widetilde{x}^0$ be defined as above. Then, there exists locally a unique continuously differentiable solution $\widetilde{x}(\lambda)$ to the PDE system stated in \eqref{eq:first_order_opt_cond2} with initial condition $\widetilde{x}(\lambda^0) = \widetilde{x}^0$ such that $\widetilde{x}(\lambda)$ is $\varepsilon$-critical with respect to $f_\lambda$ and, thus, $\varepsilon'$-Pareto critical with respect to \eqref{eq:MOP} for
    $\varepsilon' = {\varepsilon} /{\min\{\lambda_1,\dots,\lambda_p\}}>0$.
\end{theorem}

\begin{proof}
Let $c \in \R^n$ denote the vector such that $\nabla_x f_{\lambda^0}(\widetilde{x}^0) = c$. By defining $F(x,\lambda ) \coloneqq \nabla_x f_\lambda(x) - c$, it directly follows that $F(\widetilde{x}^0 ,\lambda^0) = 0$. Moreover, since $\nabla_x F(x,\lambda) = \nabla^2_x f_\lambda(x)$, it holds that $\nabla_x F(\widetilde{x}^0,\lambda^0)$ is invertible. 

Hence, as in the proof of Theorem~\ref{thm:pareto_tracing_neighborhood}, the implicit function theorem yields open neighborhoods $\U(\lambda^0), \mathcal{V}(\widetilde{x}^0)$ and a unique, continuously differentiable function $ \widetilde{x}_F \colon \mathcal{U}(\lambda^0)\to\mathcal{V}(\widetilde{x}^0)$ with $\widetilde{x}_F(\lambda^0)=\widetilde{x}^0$ such that, for all $(x,\lambda)\in  \mathcal{V}(\widetilde{x}^0)\times\U(\lambda^0)$, it holds that $\nabla_x f_\lambda(\widetilde{x}_F(\lambda))=c$ for all $\lambda\in\mathcal{U}(\lambda^0)$. Implicit differentiation of $\widetilde x(\lambda)$ with respect to $\lambda_i$, $i=1,\ldots,p-1$, yields~\eqref{eq:first_order_opt_cond2}. Uniqueness of $\widetilde{x}_F$ with $\widetilde{x}_F(\lambda^0) =\widetilde{x}^0$ can be shown analogously as in Theorem~\ref{thm:pareto_tracing_neighborhood}. Finally, we get for all $\lambda\in\mathcal{U}(\lambda^0)$,
$$
\left\|\nabla_x f_\lambda(\widetilde{x}_F(\lambda))\right\| = \|c\| = \left\|\nabla_x f_{\lambda^0}(\widetilde{x}^0)\right\| \leq \varepsilon
$$
Thus, $\widetilde{x}_F(\lambda)$ is $\varepsilon$-critical for $f_\lambda$ and $\varepsilon'$-Pareto critical with $\varepsilon' \!=\! {\varepsilon}/ {\min\{\lambda_1,\dots,\lambda_p\}}.$ \end{proof}

The local existence result in Theorem~\ref{thm:pareto_tracing_init_perturbed_local} 
can be extended to the complete projected weight space under stronger convexity assumptions.

\begin{theorem}\label{thm:pareto_tracing_init_perturbed}
Let $\lambda^0\in\proj$ and let $\widetilde{x}^0$ be defined as above, with
$c \coloneqq \nabla_x f_{\lambda^0}(\widetilde{x}^0)$.
Assume that, for every $\lambda\in\proj$, the weighted objective $f_\lambda$ is strongly convex. Define $f_\lambda^c(x):=f_\lambda(x)-x^\top c.$
Then the mapping
$$
\widetilde{x}:\proj\to\R^n,
\qquad
\widetilde{x}(\lambda)
=
\arg\min_{x\in\R^n} f_\lambda^c(x)
$$
is continuously differentiable and is the unique solution to \eqref{eq:first_order_opt_cond2} on $\proj$ with initial condition $\widetilde{x}(\lambda^0)=\widetilde{x}^0$. Moreover, $\widetilde{x}(\lambda)$ is $\varepsilon$-critical with respect to $f_\lambda$ and, consequently, $\varepsilon'$-Pareto critical for \eqref{eq:MOP}, where
$\varepsilon'={\varepsilon}/{\min\{\lambda_1,\dots,\lambda_p\}}$.

\end{theorem}

\begin{proof}
Since $f_\lambda$ is strongly convex, it follows that $f_\lambda^c$ is also strongly convex. Thus, $f_\lambda^c$ is coercive and has a unique minimum for all $\lambda\in\proj$. Moreover, its Hessian is equal to $\nabla_x^2f_\lambda(x)$, which is invertible by the strong convexity of $f_\lambda$.

We proceed as in the proof of Theorem~\ref{thm:pareto_tracing}. Let $\lambda^k\in\proj$ be arbitrary and set $\widetilde{x}^k:=\widetilde{x}(\lambda^k)=\arg\min_{x\in\R^n} f_{\lambda^k}^c(x).$ Then,
$$
\nabla_x f_{\lambda^k}^c(\widetilde{x}^k)
=
\nabla_x f_{\lambda^k}(\widetilde{x}^k)-c
=
0.
$$
Applying the implicit function theorem as in the proof of Theorem~\ref{thm:pareto_tracing_init_perturbed_local} yields a unique continuously differentiable local function $\widetilde{x}_F^k$ satisfying $\nabla_x f_\lambda^c(\widetilde{x}_F^k(\lambda))=0$. By strong convexity, $\widetilde{x}_F^k$ coincides locally with $\widetilde{x}$. Since $\lambda^k\in\proj$ was arbitrary, it follows that $\widetilde{x}\in C^1(\proj,\R^n)$. As $c$ is independent of $\lambda$, implicit differentiation yields~\eqref{eq:first_order_opt_cond2}.

Furthermore, by the definition of $c$, $\nabla_x f_{\lambda^0}^c(\widetilde{x}^0)=0$. By strong convexity, $\widetilde{x}^0$ is the unique minimizer of $f_{\lambda^0}^c$, and hence $\widetilde{x}(\lambda^0)=\widetilde{x}^0$. Uniqueness of $\widetilde{x}$ with $\widetilde{x}(\lambda^0)=\widetilde{x}^0$ follows by arguing as in Theorem~\ref{thm:pareto_tracing_neighborhood}. 

Finally, for all $\lambda\in\proj$, it holds that
$$
\left\lVert\nabla_x f_\lambda(\widetilde{x}(\lambda))\right\rVert = \lVert c\rVert
\leq\varepsilon.
$$
Hence, $\widetilde{x}(\lambda)$ is $\varepsilon$-critical with respect to $f_\lambda$ and $\varepsilon'$-Pareto critical for \eqref{eq:MOP}.
\end{proof}

The error in the solutions obtained from an approximated initial value can be bounded in terms of the initial error, as shown in the following theorem.

\begin{theorem}\label{thm:ErrorBoundInitCond} 
Let $\widetilde{x}^0$ sufficiently close to $x^0$ be defined as above, and let $x(\lambda)$ and $\widetilde{x}(\lambda)$ solve \eqref{eq:first_order_opt_cond2} locally with respect to initial values $x^0$ and $\widetilde{x}^0$, respectively. Then there exists a neighborhood $\mathcal{U}(\lambda^0) \subset \proj$ and constants $M,\mu>0$ such that
\begin{align*}
        \lVert \Tilde{x}(\lambda) - x(\lambda) \rVert \leq \frac{M}{\mu} \lVert \Tilde{x}^0 - x^0 \rVert.
    \end{align*}
\end{theorem}

\begin{proof}
In the following, let $\lVert \cdot \rVert$ denote the Euclidean norm for vectors and the corresponding induced operator norm for matrices. Using Theorems~\ref{thm:pareto_tracing_init_perturbed_local} and~\ref{thm:pareto_tracing_neighborhood}, it locally holds that $\nabla_x f_\lambda(\Tilde{x}) = c$ and $\nabla_x f_\lambda (x(\lambda)) = 0$, where $c \coloneqq \nabla_x f_{\lambda^0}(\widetilde{x}^0)$. Hence, $\nabla_x f_\lambda(\widetilde{x}(\lambda)) - \nabla_x f_\lambda(x(\lambda)) = c$. For fixed $\lambda$ in a sufficiently small 
neighborhood of $\lambda^0$, we define $g(t) = \nabla_x f_\lambda (x(\lambda) + t (\Tilde{x}(\lambda) - x(\lambda)))$ and then, we have that
\begin{align*}
    c = g(1) - g(0)  = \int_0^1 g'(t) dt &= \left( \int_0^1 \nabla_x^2 f_\lambda \left( x(\lambda) + t ( \Tilde{x}(\lambda) - x(\lambda))\right) dt  \right) (\Tilde{x}(\lambda) - x(\lambda)) \\
    & = I(\lambda) (\Tilde{x}(\lambda) - x(\lambda)).
\end{align*}

Since $\nabla_x^2 f_{\lambda^0}(x^0)$ is positive definite, the continuity of $\nabla_x^2 f_\lambda(x)$ implies that there exist $\mu,r>0$ and an open neighborhood $\mathcal{U}(\lambda^0)$ such that the smallest eigenvalue satisfies $\rho_{\min}(\nabla_x^2 f_\lambda(x))\geq\mu$  for all $\lambda\in\mathcal{U}(\lambda^0)$ and
$x\in\mathcal{B}_r(x^0)$. For $\widetilde{x}^0$ sufficiently close to $x^0$, we can choose $\mathcal{U}(\lambda^0)$ sufficiently small such that
$x(\lambda),\widetilde{x}(\lambda)\in\mathcal{B}_r(x^0)$
for all $\lambda\in\mathcal{U}(\lambda^0)$ by continuity of both solutions. By convexity, their connecting line segment is also contained in $\mathcal{B}_r(x^0)$. It directly follows that $I(\lambda)$ is uniformly positive definite on $\mathcal{U}(\lambda^0)$, i.\,e., $v^\top I(\lambda)v \geq \int_0^1 \mu \lVert v \rVert^2 dt = \mu \lVert v \rVert^2$ for all $v\in\R^n$. Therefore, $\rho_{min} (I(\lambda)) \geq \mu$, and thus $ \lVert I(\lambda)^{-1} \rVert_2 =  \rho_{max} (I(\lambda)^{-1}) = {1}/{\rho_{\min}(I(\lambda))} \leq {1}/{\mu}$ for all $\lambda \in \mathcal{U}(\lambda^0)$. Hence, we obtain 
\begin{align} \label{eq:thm_localerrorbound_eq1}
            \lVert \Tilde{x}(\lambda) - x(\lambda) \rVert = \lVert I(\lambda)^{-1}c \rVert \leq \lVert I(\lambda)^{-1} \rVert_2 \lVert c \rVert \leq \frac{1}{\mu} \lVert c \rVert.
\end{align}
It remains to bound $\lVert c\rVert$ with respect to the initial error. For $\lambda=\lambda^0$, we have $c = \nabla_x f_{\lambda^0} (\Tilde{x}^0)- \nabla_x f_{\lambda^0}(x^0)$. Analogous to the above, we define $g_0 (t) = \nabla_x f_{\lambda^0} (x^0 + t (\Tilde{x}^0 - x^0))$. Setting $I^0 \coloneqq I(\lambda^0)$, we get that $\lVert c \rVert \leq \lVert I^0 \rVert \lVert \Tilde{x}^0 - x^0 \rVert$. Since $\nabla^2_x f_{\lambda^0}$ is continuous,
we have that $\nabla^2_x f_{\lambda^0}$ is locally bounded. Thus, there exists a constant $M>0$ such that $\lVert I^0 \rVert \leq M$. Hence, 
\begin{align}\label{eq:thm_localerrorbound_eq2}
    \lVert c \rVert \leq M \lVert \Tilde{x}^0 - x^0 \rVert.
\end{align}
Using \eqref{eq:thm_localerrorbound_eq1} and \eqref{eq:thm_localerrorbound_eq2}, we finally obtain for all $\lambda\in\mathcal U(\lambda^0)$
$\lambda\in\mathcal U(\lambda^0)$
\begin{align*}
    \lVert \widetilde{x}(\lambda) - x(\lambda) \rVert  \leq \frac{1}{\mu} \lVert c \rVert \leq \frac{M}{\mu} \lVert \widetilde{x}^0 - x^0 \rVert.
\end{align*}
As all norms on $\R^n$ are equivalent, analogous estimates hold for any norm.
\end{proof}

\section{Pareto Front Tracing by Numerical Integration} \label{sec:NumIntegration}

In the previous section, we derived the theoretical framework for Pareto front tracing for \eqref{eq:MOP} with an arbitrary number of objectives based on~\eqref{eq:first_order_opt_cond2}. Since analytical solutions to this system (and thus an analytical description of the Pareto set) are generally unavailable, we present a numerical approximation scheme based on numerical integration.

\subsection{Numerical Integration and Tracing Strategy}\label{subsec:NumInteTracingStrat}

To illustrate the general idea of the proposed numerical integration approach, we first consider the tri-objective case before extending the approach to an arbitrary number of objectives.

\subsubsection{Pareto Front Tracing for Tri-objective Problems} \label{subsec:NumInteTriObjective}

\paragraph{Numerical Integration}

For the tri-objective case, the resulting PDE system is
\begin{align}
    & \partial_{\lambda_1}x(\lambda) = \, - (\nabla_x^2f_{\lambda}(x(\lambda)) )^{-1} \, (\nabla_x f_1(x(\lambda)) - \nabla_x f_3(x(\lambda)) ) \eqqcolon g_{\lambda_1}(x(\lambda), \lambda)  \label{eq:triObjPDEsystem1}\\
	& \partial_{\lambda_2}x(\lambda) = \, - (\nabla_x^2f_{\lambda}(x(\lambda)) )^{-1} \, (\nabla_x f_2(x(\lambda)) - \nabla_x f_3(x(\lambda)) ) \eqqcolon g_{\lambda_2}(x(\lambda), \lambda). \label{eq:triObjPDEsystem2}
\end{align}
Recall, that in this specific case $\lambda = (\lambda_1, \lambda_2, 1 - \lambda_1 - \lambda_2)$. A key property of this first-order PDE system is that there are no dependencies among the partial derivatives. More precisely, each equation contains exactly one partial derivative, implying that each partial derivative appears exactly once in the entire system. Consequently, each equation can be interpreted as an ordinary differential equation with respect to a single weight parameter when the remaining parameter is fixed. In the tri-objective case, this corresponds to integrating with respect to either $\lambda_1$ or $\lambda_2$. As we make use of numerical integration methods for ODEs, we briefly recall standard integration schemes for equations of the form
\begin{align} \label{eq:ODE}
    \partial_\lambda x(\lambda) = g(x(\lambda),\lambda),
\end{align}
which, in the context of Pareto front tracing, corresponds exactly to the differential equation obtained in the bi-objective case.

One approach is the Runge--Kutta method~\citep{hair:solvODE:1993}, which computes approximations by
\begin{align}
    &x(\lambda^j+h)
    \approx x^{j+1}
    \coloneqq x^j+h\,\sum_{r=1}^{s} b_r k_r,
    \label{eq:runge_kutta}\\
   & k_r
    \coloneqq
    g\left(
        x^j+h\,\sum_{\ell=1}^{r-1}a_{r\ell}k_\ell,
        \lambda^j+c_rh
    \right),
    \qquad r=1,\dots,s, 
\end{align}
where $a_{r\ell}$, $b_r$, and $c_r$ denote the coefficients of the
RK method. The RK scheme can be applied straightforwardly to integrate the two equations of the PDE system with respect to $\lambda_1$ and $\lambda_2$, respectively:

\begin{align}
    \Phi_{\lambda_i}^h(x^{j,k})
    \coloneqq
    x^{j,k}+h\sum_{r=1}^s b_r k_{r,i}^{j,k}, \quad\text{where}
    \label{eq:RKPhi}
\end{align}
\begin{align}
k_{r,i}^{j,k}\coloneqq g_{\lambda_i}\left( x^{j,k}+h\sum_{\ell=1}^{r-1}a_{r\ell}k_{\ell,i}^{j,k}, \lambda^{j,k}+c_r h e_i\right),
\quad r=1,\dots,s, \quad i=1,2.
\end{align}

The resulting RK updates are then given by
\begin{align}
    x(\lambda_1^{j}+h,\lambda_2^k)
    \approx x^{j+1,k}
    &\coloneqq \Phi_{\lambda_1}^h(x^{j,k}),
    \label{eq:RK1}\\
    x(\lambda_1^{j},\lambda_2^k+h)
    \approx x^{j,k+1}
    &\coloneqq \Phi_{\lambda_2}^h(x^{j,k}).
    \label{eq:RK2}
\end{align}

\paragraph{Tracing Strategy}

In practice, we usually want each objective function to be weighted with a minimum value of $\delta > 0$ 
to bound the trade-off between function values. Therefore, we replace the weight space $\Lambda$ by a closed subset $\Lambda_\delta = \left\{ \lambda \in \mathbb{R}^p : \lambda_i \geq \delta \text{ for } i=1, \dots, p,  \, \sum_{i=1}^{p} \lambda_i = 1 \right\}$, where $0 < \delta < {1}/{p}$. We refer to $\Lambda_\delta$ as the \textit{$\delta$-weight space}. Analogously, we denote its projection onto $\R^{p-1}$ by $\proj_\delta$, which we also refer to as the $\delta$-weight space, since each $\lambda\in\proj_\delta$ uniquely corresponds to a weight vector in $\Lambda_\delta$. In the tri-objective case, where $\lambda_3 = 1-\lambda_1-\lambda_2$, we have $\proj_\delta = \{(\lambda_1,\lambda_2)\in\R^2 \colon \lambda_1,\lambda_2\geq\delta,\; \lambda_1+\lambda_2\leq1-\delta\}$. Restricting the weight space to $\Lambda_\delta$ provides a compact computational domain that can be discretized by a finite grid.

We define the grid $\mathcal{G}$ on $[\delta,1-2\delta]^2$ by $ \mathcal{G} = \{\lambda_1^0,\dots,\lambda_1^t\} \times \{\lambda_2^0,\dots,\lambda_2^t\},$ where $\lambda_i^j=\delta+jh$ for $j=0,\dots,t$ and $i=1,2$, and $h>0$ denotes the step size. Note that $\lambda_i^t = \delta + th$ is well-defined if  $(1-3\delta)/h = t \in \N$. In the following, we set $h = \delta$ to simplify the presentation. 
Finally, let $\mathcal{G}_{\Lambda_\delta}=\mathcal{G}\cap\proj_\delta$ denote the set of all grid points $\lambda^{j,k}$ at which approximations $x^{j,k}$ are computed.

The proposed tracing strategy 
consists of three stages. Starting from an initial $\varepsilon$-critical solution $x^{0,0}$ at $\lambda^{0,0}=(\delta,\delta)$, we first integrate along both coordinate directions using \eqref{eq:RK1} and \eqref{eq:RK2}, respectively, to obtain approximations $x^{j,0}$ and $x^{0,k}$, for $j,k=1,\dots,t$. We refer to these approximations, together with the initial solution $x^{0,0}$, as the boundary $\mathcal{B}(\mathcal{G}_{\Lambda_\delta},x^{0,0})$ of the tracing strategy with respect to $x^{0,0}$.

In the second stage, the generated approximations in $\mathcal{B}(\mathcal{G}_{\Lambda_\delta},x^{0,0})$ serve as initial values for integrations in the respective other coordinate direction until the entire grid $\mathcal{G}_{\Lambda_\delta}$ is covered. This yields two approximations $x_{\lambda_1}^{j,k}$ and $x_{\lambda_2}^{j,k}$ at interior grid points, obtained by integration in $\lambda_1$- and $\lambda_2$-direction, respectively. In the third stage, when both directional approximations are available, their arithmetic mean is computed.

Note that, although averaging could be performed directly during the second stage for fixed step sizes, we formulate it as a separate third stage, as the practical implementation makes use of adaptive step-size solvers, see Section~\ref{subsec:expsetup}.

In general, the tracing procedure can also be initialized at a different weight vector. However, in contrast to the bi-objective case, choosing an initial weight vector in the interior of the grid does not provide an analogous reduction in the required integration path lengths. A different initial weight may nevertheless be preferable if the weighted Hessian is poorly conditioned at $\lambda^{0,0}=(\delta,\delta)$.

\begin{figure}[t]
	\begin{minipage}[t]{0.32\textwidth}
		\centering
			\begin{tikzpicture}[scale=0.345]

				\draw[line width=0.6pt,dotted,color=black!90] (0,0) -- (9,0);
				\draw[line width=0.6pt,dotted,color=black!90] (0,0) -- (0,9);
				\draw[line width=0.6pt,dotted,color=black!90] (9,0) -- (0,9);
				
				\foreach \a/\b in {1.5/7.5, 3/6, 4.5/4.5, 6/3, 7.5/1.5} {
					\draw[line width=0.6pt,dotted,color=black!90] (0,\a) -- (\b,\a);
					\draw[line width=0.6pt,dotted,color=black!90] (\a,0) -- (\a,\b);
				}

				\filldraw[color=mediumblue] (0,0) circle (2pt) node[below left] {\tiny $\color{mediumblue}\lambda^{0,0}$};

				\draw[line width=0.8pt,->,color=mediumblue] (0,0) -- (8.8,0);

				\draw[line width=0.8pt,->,color=mediumblue] (0,0) -- (0,8.8);

				\filldraw[color=mediumblue] (9,0) circle (2pt) node[below] {\tiny $\color{mediumblue}\lambda^{t,0}$};
				\filldraw[color=mediumblue] (0,9) circle (2pt) node[left] {\tiny $\color{mediumblue}\lambda^{0,t}$};
				
			\end{tikzpicture}

		\subcaption{Stage 1: Computation of the boundary
		$\mathcal{B}(\mathcal{G}_{\Lambda_\delta},x^{0,0})$ (blue).}
		\label{fig:TracingStrat1}
	\end{minipage}
	\hfill
	\begin{minipage}[t]{0.32\textwidth}
		\centering
			\begin{tikzpicture}[scale=0.345]

				\draw[line width=0.6pt,dotted,color=black!90] (0,0) -- (9,0);
				\draw[line width=0.6pt,dotted,color=black!90] (0,0) -- (0,9);
				\draw[line width=0.6pt,dotted,color=black!90] (9,0) -- (0,9);
				
				\foreach \a/\b in {1.5/7.5, 3/6, 4.5/4.5, 6/3, 7.5/1.5} {
					\draw[line width=0.6pt,dotted,color=black!90] (0,\a) -- (\b,\a);
					\draw[line width=0.6pt,dotted,color=black!90] (\a,0) -- (\a,\b);
				}

				\filldraw[color=mediumblue] (0,0) circle (2pt) node[below left] {\tiny $\color{mediumblue}\lambda^{0,0}$};

				\foreach \y/\name in {1.5/1, 3/2, 4.5/3, 6/4} {
					\filldraw[color=magenta] (0,\y) circle (2pt) node[left] {\tiny $\color{magenta}\lambda^{0,\name}$};
				}

                \filldraw[color=magenta] (0,7.5) circle (2pt);
                \filldraw[color=magenta] (0,7.75) node[left] {\tiny $\color{magenta} \vdots $};

				\foreach \y/\xend in {1.5/7.3, 3/5.8, 4.5/4.3, 6/2.8, 7.5/1.3} {
					\draw[line width=0.8pt,->,color=magenta] (0,\y) -- (\xend,\y);
				}

				\foreach \x/\name in {1.5/1, 3/2, 4.5/3, 6/4} {
					\filldraw[color=forestgreen(web)] (\x,0) circle (2pt) node[below] {\tiny $\color{forestgreen(web)}\lambda^{\name,0}$};
				}

                \filldraw[color=forestgreen(web)] (7.5,0) circle (2pt);
                \filldraw[color=forestgreen(web)] (7.6,0) node[below] {\tiny $\color{forestgreen(web)} \dots $};

				\foreach \x/\yend in {1.5/7.3, 3/5.8, 4.5/4.3, 6/2.8, 7.5/1.3} {
					\draw[line width=0.8pt,->,color=forestgreen(web)] (\x,0) -- (\x,\yend);
				}

				\filldraw[color=mediumblue] (9,0) circle (2pt) node[below] {\tiny $\color{mediumblue}\lambda^{t,0}$};
				\filldraw[color=mediumblue] (0,9) circle (2pt) node[left] {\tiny $\color{mediumblue}\lambda^{0,t}$};
				
			\end{tikzpicture}
		
        \subcaption{Stage 2: Integration from $\mathcal{B}(\mathcal{G}_{\Lambda_\delta},x^{0,0})$ in
		$\lambda_1$- (magenta) and $\lambda_2$-direction (green).}
		\label{fig:TracingStrat2}
	\end{minipage}
	\hfill
	\begin{minipage}[t]{0.32\textwidth}
		\centering
			\begin{tikzpicture}[scale=0.345]

				\draw[line width=0.6,-,dotted, color = black!90] (0,0) -- (9,0) ;
				\draw[line width=0.6,-,dotted, color = black!90] (0,0) -- (0,9) ;
				\draw[line width=0.6,-,dotted, color = black!90] (9,0) -- (0,9) ;
				
				\foreach \a/\b in { 1.5/7.5, 3/6, 4.5/4.5, 6/3}
				\draw[line width=0.6,-,dotted, color = black!90] (\a,0) -- (\a,\b) ;
				\draw[line width=0.6,-,dotted, color = black!90] (7.5,0.2) -- (7.5,1.5) ;
				\foreach \a/\b in { 1.5/7.5, 3/6, 4.5/4.5, 6/3}
				\draw[line width=0.6,-,dotted, color = black!90] (0,\a) -- (\b,\a) ;
				\draw[line width=0.6,-,dotted, color = black!90] (0.2,7.5) -- (1.5,7.5) ;

				\filldraw[color = mediumblue] (0,0) circle (2pt) node[below left] {\tiny $\color{mediumblue} \lambda^{0,0}$};

                \filldraw[color=magenta] (0,7.75) node[left] {\tiny $\color{magenta} \vdots $};
                \filldraw[color=forestgreen(web)] (7.6,0) node[below] {\tiny $\color{orange} \dots $};
				
				\filldraw[color = magenta] (0,9) circle (2pt) node[left] {\tiny $\color{magenta} \lambda^{0,t}$};
				\filldraw[color = forestgreen(web)] (9,0) circle (2pt) node[below] {\tiny $\color{forestgreen(web)} \lambda^{t,0}$};
				
				\foreach \x/\y/\name in {0/1.5/1, 0/3/2, 0/4.5/3, 0/6/4}
				\filldraw[color = magenta] (\x,\y) circle (2pt) node[left] {\tiny $\color{magenta} \lambda^{0,\name}$};
				\filldraw[color = magenta] (0,7.5) circle (2pt);
				
				\foreach \x/\y/\name in {1.5/0/1, 3/0/2, 4.5/0/3, 6/0/4}
				\filldraw[color = forestgreen(web)] (\x,\y) circle (2pt) node[below] {\tiny $\color{forestgreen(web)} \lambda^{\name,0}$};
				\filldraw[color = forestgreen(web)] (7.5,0) circle (2pt);
				
				\foreach \x/\y in {1.5/1.5, 1.5/3, 1.5/4.5, 1.5/6, 1.5/7.5, 3/3, 4.5/3, 6/3, 3/1.5, 3/4.5, 3/6, 4.5/1.5, 4.5/4.5, 6/1.5, 7.5/1.5}
				\filldraw[color = orange] (\x,\y) circle (2pt);
				
			\end{tikzpicture}
		
        \subcaption{Stage 3: Averaged approximations at the interior grid
		points (orange).}
        \label{fig:TracingStrat3}
	\end{minipage}
	\caption{Illustration of the three stages of the tracing strategy
	for tri-objective problems.}
	\label{fig:TracingStrat}
\end{figure}
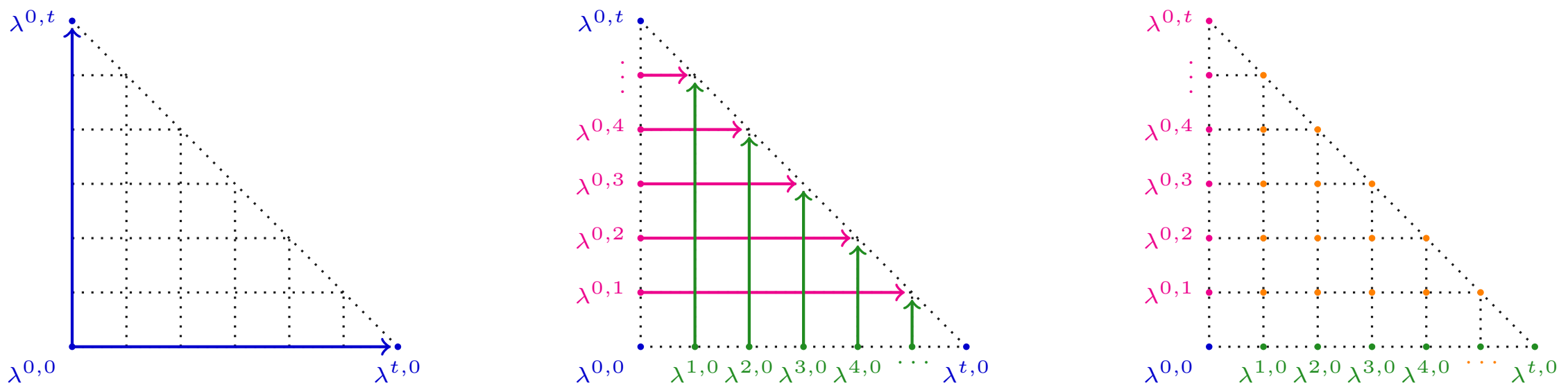

Algorithm~\ref{alg:tracing_tri_obj} summarizes the proposed tracing strategy for tri-objective problems. The three stages are illustrated in Figure~\ref{fig:TracingStrat}.

\begin{algorithm}[tbp]
\caption{Pareto front tracing for tri-objective problems}
\label{alg:tracing_tri_obj}
\begin{algorithmic}[1]

\REQUIRE Parameter $\delta$, initial value $x^{0,0}$ for $\lambda^{0,0}=(\delta,\delta)$, RK updates $\Phi_{\lambda_i}^h$, $i=1,2$
\ENSURE Approximated values $x^{j,k}$ of Pareto critical solutions

\STATE $h=\delta$, $t=\frac{1}{\delta}-3$

\COMMENT{Stage 1: Compute boundary approximations}

\FOR{$j=0,\dots,t-1$}
    \STATE $x^{j+1,0}=\Phi_{\lambda_1}^h(x^{j,0})$, \quad
           $x^{0,j+1}=\Phi_{\lambda_2}^h(x^{0,j})$
\ENDFOR

\COMMENT{Stage 2: Compute directional approximations}
\STATE $x_{\lambda_1}^{0,k}=x^{0,k},\quad
        x_{\lambda_2}^{j,0}=x^{j,0}$,
        \quad $j,k=1,\dots,t-1$
\FOR{$k=1,\dots,t-1$}
    \FOR{$j=1,\dots,t-k$}
        \STATE $x_{\lambda_1}^{j,k}
        =\Phi_{\lambda_1}^h(x_{\lambda_1}^{j-1,k})$
        \STATE $x_{\lambda_2}^{j,k}
        =\Phi_{\lambda_2}^h(x_{\lambda_2}^{j,k-1})$
    \ENDFOR
\ENDFOR

\COMMENT{Stage 3: Average interior approximations}
\FOR{$k=1,\dots,t-1$}
    \FOR{$j=1,\dots,t-k$}
        \STATE $x^{j,k}=\frac{1}{2}
        \left(x_{\lambda_1}^{j,k}+x_{\lambda_2}^{j,k}\right)$
    \ENDFOR
\ENDFOR

\end{algorithmic}
\end{algorithm}

\subsubsection{Extension to an Arbitrary Number of Objectives}

The tracing strategy introduced for tri-objective problems can be generalized to an arbitrary number of $p$ objective functions by recursion. The normalized weight
vector is determined by $d=p-1$ independent weight parameters $\lambda_1,\dots,\lambda_d$. Starting from $\lambda^0=(\delta,\dots,\delta)$ in the projected $\delta$-weight space, we define the tracing boundary as the union of the $d$ facets obtained by fixing
one of the independent weight parameters at $\delta$.

Each of these facets has dimension $d-1$ and can therefore be traced by
recursively applying the same procedure to the remaining independent weight
parameters. More precisely, for every $i\in\{1,\dots,d\}$, the facet $\mathcal{F}_i \coloneqq \{\lambda\in\mathcal{G}_{\Lambda_\delta} \colon \lambda_i=\delta\}
$ is computed by applying the tracing strategy with $\lambda_i$ fixed and
$\lambda_j$, $j\neq i$, as the remaining free parameters. The recursion
terminates in the bi-objective case, where only one independent weight
parameter remains and the corresponding ODE is integrated directly.

Once approximations have been obtained on all facets
$\mathcal{F}_i$, each of them serves as a collection of initial values for
integrations in the corresponding $\lambda_i$-direction. In this way, an
interior grid point can be approximated independently from each of the $d$
coordinate directions. The final approximation is defined as the arithmetic
mean of these $d$ approximations.

\begin{remark}
The proposed tracing strategy yields approximations of Pareto optimal solutions for a predefined grid of the $\delta$-weight space. Approximations for weight vectors in a grid cell can be obtained using interpolation or convex combinations.
\end{remark}

\subsection{Error Analysis}

For an error analysis of the proposed tracing strategy, we consider the individual ODEs arising from the PDE system \eqref{eq:first_order_opt_cond2} by fixing all weight parameters except $\lambda_i$, as discussed in Section~\ref{subsec:NumInteTracingStrat}. 

For fixed values $\bar\lambda_b$, $b\neq i$, let $\bar\lambda_{-i}$ denote the vector of these fixed parameters and let $I_i(\bar\lambda_{-i})=[a_i(\bar\lambda_{-i}),b_i(\bar\lambda_{-i})]$ be a compact interval such that the corresponding weight vectors $\lambda(t)=(\bar\lambda_1,\ldots,\bar\lambda_{i-1},t,\bar\lambda_{i+1},\ldots,\bar\lambda_{p-1})$ remain in $U(\lambda^0)\subset\proj$, on which the solution $x$ exists, for all $t\in I_i(\bar\lambda_{-i})$. Setting $y_{i,\bar\lambda_{-i}}(t) \coloneqq x(\lambda(t))$ and defining $g_i(z,t) \coloneqq g_{\lambda_i}(z, \lambda(t))$, we get by using the chain rule
\begin{align}\label{eq:partialDerivative_ODE_form}
\partial_t y_{i,\bar\lambda_{-i}} = \partial_{\lambda_i}x(\lambda(t))= g_{\lambda_i}(y_{i,\bar\lambda_{-i}},\lambda(t)) = g_i (y_{i,\bar\lambda_{-i}},t) , \qquad t \in I_i(\bar\lambda_{-i}),
\end{align}
which yields the corresponding ODEs in standard form.

Based on these, we consider an integration path as introduced in Section~\ref{subsec:NumInteTriObjective}, which involves successively solving several individual ODEs on compact subintervals $[t_{i_0},t_{i_{N}}]\subseteq I_i(\bar\lambda_{-i})$. To analyze the error propagation along the path, we first derive local error bounds that hold uniformly for all coordinate directions. This requires a suitable compact set $\mathcal C\subset V(x^0)\times U(\lambda^0)$ on which all $g_{\lambda_i}$, $i=1,\ldots,p-1$, are well-defined. Such a set can be chosen as a sufficiently small compact neighborhood of the exact solution along the entire integration path, i.\,e., of the set $\{(x(\lambda),\lambda) \mid \lambda\in C_\Lambda\}$, where $C_\Lambda \coloneqq \bigcup_{i=1}^{p-1}\{\lambda(t) \mid t \in [t_{i_0},t_{i_N}]\}$. Note that, for notational convenience, we restrict the error analysis to forward integration. Analogous error estimates for backward integration can be derived similarly.

\paragraph{Perturbations of Right-hand Sides of the PDE System}

While the effect of an input error caused by approximated initial values is examined in Section~\ref{subsec:InitCond}, the effects of perturbations in the right-hand sides of the PDE system~\eqref{eq:first_order_opt_cond2} have not yet been discussed. In many practical applications, $g_{\lambda_i}$, $i \in \{1, \dots, p-1\}$ is not given explicitly, but has to be approximated by a sequence $(g_{\lambda_i, \ell})_{\ell \in \N}$, which we assume to converge uniformly to $g_{\lambda_i}$ on $\mathcal{C}$, i.\,e., $\sup_{(z,\lambda)\in\mathcal C} \left\| g_{\lambda_i}(z,\lambda)-g_{\lambda_i,\ell}(z,\lambda)\right\| \to 0$  as $\ell \to \infty$. For any $\varepsilon_i > 0$, we choose $\ell^*$ sufficiently large such that $ \sup_{(z,\lambda)\in\mathcal C} \left\|g_{\lambda_i}(z,\lambda) - g_{\lambda_i,\ell^*}(z,\lambda)\right\|   \leq \varepsilon_i $ for all $i=1,\ldots,p-1$. Define $\widetilde {g}_{\lambda_i} \coloneqq g_{\lambda_i,\ell^*}$ and $\widetilde{\varepsilon} \coloneqq \max_{i=1,\ldots,p-1}\varepsilon_i$. The following estimate holds for each of the individual ODEs considered along the integration path.

\begin{lemma}[cf. {\citep[Theorem I.12.1]{agar:introode:2008}}]\label{lem:perturbedRHS_local}
Let $[t_{i_0}, t_{i_N}] \subseteq I_i(\bar\lambda_{-i})$. Assume that there exists $L>0$ such that
$$
\|g_{\lambda_i}(z^1,\lambda)-g_{\lambda_i}(z^2,\lambda)\|
\leq L\|z^1-z^2\|
$$
for all $i=1,\ldots,p-1$ and $(z^1,\lambda),(z^2,\lambda)\in\mathcal{C}$. Let $\widetilde{g}_{\lambda_i}$ be as defined above and suppose that $\widetilde{g}_{\lambda_i}$ is locally Lipschitz continuous on $\operatorname{int}(\mathcal{C})$. For sufficiently small $\widetilde\varepsilon$ and an initial value $\widetilde y_{i,\bar\lambda_{-i}}(t_{i_0})$ sufficiently close to $y_{i,\bar\lambda_{-i}}(t_{i_0})$, the solution $\widetilde y_{i,\bar\lambda_{-i}}$ of~\eqref{eq:partialDerivative_ODE_form} with perturbed right-hand side $\widetilde{g}_{\lambda_i}$ exists on $[t_{i_0},t_{i_N}]$ and satisfies 
\begin{align*}
\left\| y_{i,\bar\lambda_{-i}}(t_{i_{N}}) \! -\! \widetilde y_{i,\bar\lambda_{-i}}(t_{i_{N}})\right\| 
\leq  \, e^{L(t_{i_{N}}\!-\!t_{i_0})} \left\| y_{i,\bar\lambda_{-i}}(t_{i_0}) -\widetilde y_{i,\bar\lambda_{-i}}(t_{i_0}) \right\| \! + \! \frac{\widetilde\varepsilon}{L} \left(e^{L(t_{i_{N}} \! -\! t_{i_0})}-1\right)
\end{align*}
\end{lemma}

\begin{proof}

Analogous to the proof of \cite[Proposition~2.5]{bolt:trac:2021}.
\end{proof}

Since the value at the endpoint of each coordinate integration serves as the initial value for the subsequent ODE, we can recursively apply Lemma~\ref{lem:perturbedRHS_local} to obtain a global error bound for a single integration path controlled by the initial error and the perturbations of the right-hand sides of the PDE system.

\begin{theorem}\label{thm:perturbedRHS_global}
Assume that the assumptions of Lemma~\ref{lem:perturbedRHS_local} hold for each of the individual ODEs of the form~\eqref{eq:partialDerivative_ODE_form} in the integration path with Lipschitz constant $L>0$. For sufficiently small $\widetilde\varepsilon$ and an
initial value $\widetilde x^0$ sufficiently close to $x^0$, the approximated solution exists along the entire integration path, and the approximation $\widetilde x^N$
at $\lambda^N$ satisfies
$$ \|x(\lambda^N)-\widetilde x^N\| \leq e^{L\sum_{i=1}^{p-1}(t_{i_N}-t_{i_0})} \|x^0 - \widetilde{x}^0 \| +\frac{\widetilde\varepsilon}{L} \left(e^{L\sum_{i=1}^{p-1}(t_{i_N}-t_{i_0})}-1 \right).
$$
\end{theorem}

\begin{proof}
Since the integration path consists of finitely many coordinate integrations, repeated application of Lemma~\ref{lem:perturbedRHS_local} yields existence.

Let $E_0 \coloneqq \|x^0-\widetilde{x}^0\|$ denote the initial error and $E_i$ the error at the endpoint of the $i$-th coordinate integration, $i=1,\ldots,p-1$. W.l.o.g., assume that the integration path from $\lambda^0$ to $\lambda^N$ considers the coordinate directions in ascending order and entire intervals $[t_{i_0},t_{i_N}]=[\lambda_i^0,\lambda_i^N]$ at once. Using Lemma~\ref{lem:perturbedRHS_local} and $t_{i_0} = t_{(i-1)_N}$, we get
\begin{align*}
E_i \!=\! \left\|y_{i,\bar\lambda_{-\!i}}(t_{i_N}\!) \!-\! \widetilde y_{i,\bar\lambda_{-\!i}}(t_{i_N}\!)\right\| \! & \leq   e^{L(t_{i_N}\!- t_{i_0})} \! \left\|y_{i,\bar\lambda_{-\!i}}(t_{i_0}) \!-\!\widetilde y_{i,\bar\lambda_{-\!i}}(t_{i_0})\right\|\!+\!\frac{\widetilde\varepsilon}{L} \! \left( \! e^{L(t_{i_N}-t_{i_0})}\!-\!1 \!\right)\\
&=e^{L(t_{i_N}-t_{i_0})}E_{i-1} +\frac{\widetilde\varepsilon}{L} \left(e^{L(t_{i_N}\!-t_{i_0})}-1\right)
\end{align*}
for $i=1,\ldots,p-1$. Proceeding recursively, we obtain
\begin{align*}
E_{p-1}
&\leq e^{L\sum_{i=1}^{p-1}(t_{i_N}-t_{i_0})}E_0 +\frac{\widetilde\varepsilon}{L} \left( e^{L\sum_{i=1}^{p-1}(t_{i_N}-t_{i_0})}-1 \right) .
\end{align*}
\end{proof}

\paragraph{Numerical Integration Error Analysis} \label{subsubsec:RK_Error_Ana}
Next, we examine the error induced by numerical integration using RK methods \eqref{eq:RK1} and \eqref{eq:RK2}. We assume exact initial values and no perturbations of the right-hand sides. Existing error estimates 
rely on the regularity of the exact solution, so we recall the following result.

\begin{lemma}\label{lem:regularity} 
Assume that $f_i\in C^{q+2}$ for all $i=1,\ldots,p$, $q\in\N_0$. Let $I_i(\bar\lambda_{-i})=[a_i(\bar\lambda_{-i}),b_i(\bar\lambda_{-i})]$ be defined as above. Then $y_{i,\bar\lambda_{-i}}$ is $q+1$ times continuously differentiable on $I_i(\bar\lambda_{-i})$, with bounded $(q+1)$-st derivative.
\end{lemma}

\begin{proof}
Since $f_i\in C^{q+2}$ for all $i=1,\ldots,p$ and the inversion of an invertible matrix $A$ is $C^{\infty}$ in a neighborhood of $A$, it follows that $g_{\lambda_i} \in C^q$. Using continuity of $x$ and $\lambda(t) \in C^{\infty}$, $g_{\lambda_i}(x(\lambda(t)),\lambda(t))$ is continuous as well. Hence, by \eqref{eq:partialDerivative_ODE_form}, $y_{i,\bar\lambda_{-i}}\in C^1$. Repeated application of the chain rule and induction yield $y_{i,\bar\lambda_{-i}}\in C^{q+1}$. As $I_i(\bar\lambda_{-i})$ is compact, continuity of the $(q+1)$-st derivative implies its boundedness on $I_i(\bar\lambda_{-i})$.
\end{proof}

As discussed in \cite[Lemma~2.4]{bolt:trac:2021} the regularity result also holds under weaker differentiability assumptions. Using this, we assume $f_i\in C^{q+2}$ in the following.

For sufficiently smooth problems of the form~\eqref{eq:partialDerivative_ODE_form}, the accuracy of RK methods depends on their respective orders. A RK method is said to be of order $q$ if there exists a constant $K$ independent of the step size $h$ such that
\begin{align*}
    \lVert y_{i,\bar\lambda_{-i}}(t_0+h)-y_{i,\bar\lambda_{-i}}^1\rVert \leq K h^{q+1},
\end{align*}
where $y_{i,\bar\lambda_{-i}}^1$ is obtained by applying one RK step~\eqref{eq:runge_kutta} starting from 
$y_{i,\bar\lambda_{-i}}^0\!=y_{i,\bar\lambda_{-i}}(t_0)$.

\begin{theorem}[cf. {\citep[Theorem II.3.1]{hair:solvODE:1993}}]\label{thm:RK_ODE_localError}
If the RK method as presented in \eqref{eq:runge_kutta} is of order $q \geq 1$ and if $g_i \in C^q$ for $i \in \{1, \dots, p-1\}$, then it holds that
\begin{align*}
& \lVert y_{i,\bar\lambda_{-i}}(t_0+h)-y_{i,\bar\lambda_{-i}}^1\rVert \\
\leq  & \ h^{q+1}\Bigg(\frac{1}{(q+1)!} \max_{\theta\in[0,1]} \lVert y_{i,\bar\lambda_{-i}}^{(q+1)}(t_0+\theta h)\rVert +\frac{1}{q!}\sum_{r=1}^{s}|b_r|
\max_{\theta\in[0,1]} \lVert k_r^{(q)}(\theta h)\rVert \Bigg),
\end{align*}
where $k_r^{(q)}$ denotes the $q$-th derivative of the $r$-th RK stage with respect to $h$.
\end{theorem}

By Lemma~\ref{lem:regularity}, the assumptions of Theorem~\ref{thm:RK_ODE_localError} are satisfied if $f_i\in C^{q+2}, i=1,\ldots,p$.
Indeed, the estimate holds for each step of the RK method~\citep{hair:solvODE:1993}. To derive a global error bound for the proposed tracing strategy, we show that the local errors can be bounded uniformly for all coordinate directions along an integration path. 

\begin{lemma}\label{lem:RK_uniformLocalError}
Let the RK method in \eqref{eq:runge_kutta} be of order $q\geq 1$ and $g_{\lambda_i}\in C^q$ for all $i=1,\ldots,p-1$. Consider a sequence of RK steps with fixed step size $h>0$, and let $i_\nu\in\{1,\ldots,p-1\}$ denote the coordinate integration direction which is used to compute $x^{\nu+1}$ from $x^\nu$, $\nu=0,\ldots,N-1$. For each step $\nu$, set $t_\nu \coloneqq \lambda_{i_\nu}^{\nu}$ and let $\hat y_{i_\nu,\bar\lambda_{-i_\nu}}$ denote the exact solution of the corresponding ODE with initial condition $\hat y_{i_\nu,\bar\lambda_{-i_\nu}}(t_\nu)=x^{\nu}$. Then there exists a constant $K>0$ such that 
$$
\left\| e_{\nu+1} \right\| = \left\| \hat y_{i_\nu,\bar\lambda_{-i_\nu}} (t_\nu+h)-x^{\nu+1}\right\| \leq Kh^{q+1}, 
$$ 
for all $\nu=0,\ldots,N-1$, where $h$ is sufficiently small such that the
approximations, the corresponding exact local solutions, and all RK stage evaluations remain in $\mathcal{C}$.
\end{lemma}
\begin{proof}
First, we apply Theorem~\ref{thm:RK_ODE_localError} to each RK step and obtain
\begin{align*}
\|e_{\nu+1}\| \leq h^{q+1}\Bigg( &\frac{1}{(q+1)!} \max_{\theta\in[0,1]} \left\| \hat y_{i_\nu,\bar\lambda_{-i_\nu}}^{(q+1)}(t_\nu+\theta h) \right\| +\frac{1}{q!}\sum_{r=1}^s |b_r| \max_{\theta\in[0,1]} \|k_{r,i_\nu}^{(q)}(\theta h)\| \Bigg).
\end{align*}
For $i=1,\ldots,p-1$, we define $G_{i,1}(z,\lambda)\coloneqq g_{\lambda_i}(z,\lambda)$
and recursively set $G_{i,m+1}(z,\lambda) \coloneqq D_zG_{i,m}(z,\lambda)\,g_{\lambda_i}(z,\lambda) +\partial_{\lambda_i}G_{i,m}(z,\lambda)$ for $m=1,\ldots,q$. Then, by repeatedly using the chain rule, we get $\hat y_{i,\bar\lambda_{-i}}^{(q+1)}(t) = G_{i,q+1} (\hat y_{i,\bar\lambda_{-i}}(t),\lambda(t))$. Since $g_{\lambda_i}\in C^q$, we can easily see that $G_{i,q+1}$, $i = 1, \dots, p-1$ is continuous as well. Hence, as $\mathcal{C}$ is compact, we get
$$
\left\|\hat y_{i_\nu,\bar\lambda_{-i_\nu}}^{(q+1)}(t)\right\| = \left\|G_{i_\nu,q+1} \left(\hat y_{i_\nu,\bar\lambda_{-i_\nu}}(t), \lambda(t)\right)\right\| \leq \max_{i=1,\ldots,p-1} \max_{(z,\lambda)\in \mathcal{C}} \|G_{i,q+1}(z,\lambda)\| \eqqcolon K_1
$$
for all $t\in[t_\nu,t_\nu+h]$ and $\nu=0,\ldots,N-1$. Since, moreover, the RK stages are defined recursively by $g_{\lambda_i}$, their $q$-th derivatives can also be expressed by derivatives of $g_{\lambda_i}$ up to order $q$. Analogous to the above, by continuity and compactness, we obtain
$$
\max_{\theta\in[0,1]} \lVert k_{r,i_\nu}^{(q)}(\theta h)\|
\leq
 K_2
$$
for all $\nu=0,\ldots,N-1$ and $r=1,\ldots,s$. Finally, combining both estimates yields
$$
\|e_{\nu+1}\| \leq h^{q+1} \left( \frac{K_1}{(q+1)!} +\frac{K_2}{q!} \sum_{r=1}^s |b_r| \right)\eqqcolon Kh^{q+1}.
$$
Thus, there exists a constant $K$ that yields a local error bound for all RK steps.
\end{proof}

Let $N$ denote the total number of RK steps along an integration path and $x^N$ the resulting approximation at $\lambda^N$. We define the global error by $E\coloneqq x(\lambda^N)-x^N$, which can be bounded by using Lemma~\ref{lem:RK_uniformLocalError}.

\begin{theorem}
\label{thm:RK_globalError}
Assume that the conditions of Lemma~\ref{lem:RK_uniformLocalError} are satisfied with local error constant $K>0$. 
Further, assume that there exists $L>0$ such that $ \|D_z g_{\lambda_i}(z,\lambda)\|\leq L$ for all $(z,\lambda)\in \mathcal{C}$ and $i=1,\ldots,p-1$. Then the global error can be bounded by
\begin{align}\label{eq:RK_GlobalErrorBound}
    \|E\|\leq h^q\frac{K}{L} \left(e^{LNh}-1\right),
\end{align}
where $h$ is sufficiently small for the numerical approximations to remain in $\mathcal{C}$.
\end{theorem}

\begin{proof}
Given the local error bound from Lemma~\ref{lem:RK_uniformLocalError}, we proceed as in the proof of \cite[Theorem II.3.4]{hair:solvODE:1993}. Therefore, let $E_\nu$ denote the error at $\lambda^N$ between the exact solutions with initial conditions $(\lambda^\nu,x^\nu)$ and $(\lambda^{\nu-1},x^{\nu-1})$, where $x^\nu$ and $x^{\nu-1}$ are the approximations obtained after $\nu$ and $\nu-1$ RK steps, respectively. The same Lipschitz constant $L$ and local error bound apply to all integration directions, so we obtain 
\[
\|E_\nu\| \leq e^{L(N-\nu)h}\|e_\nu\|.
\]
As in the proof of \cite[Theorem II.3.4]{hair:solvODE:1993}, we obtain
\begin{align*}
\|E\| \leq \sum_{\nu=1}^{N}\|E_\nu\|  & \leq Kh^{q+1} \sum_{\nu =1}^{N}e^{L(N-\nu)h}  = Kh^q\left( h\sum_{\nu=1}^{N}e^{L(N-\nu)h}\right)\\
&   \leq Kh^q\int_0^{Nh}e^{L(Nh-u)}\,du =h^q\frac{K}{L}(e^{LNh}-1),
\end{align*}
which gives the global error bound in~\eqref{eq:RK_GlobalErrorBound}.
\end{proof}

\begin{remark}
By the triangle inequality, the error estimates in Theorems~\ref{thm:perturbedRHS_global} and~\ref{thm:RK_globalError} also hold for the averaged approximations in Stage~3 of the tracing strategy, as the respective error bounds hold uniformly for the considered integration paths.
\end{remark}

\begin{remark}
The total tracing error can be bounded by 
the perturbation and numerical integration error using the triangle inequality, provided the assumptions of Theorem~\ref{thm:RK_globalError} also hold for the perturbed right hand sides considered in Theorem~\ref{thm:perturbedRHS_global}.
\end{remark}

\paragraph{Error of Convex Combinations}

Let $Q\subset U(\lambda^0)$ be a fixed compact set and $Q(h)\subset Q$ be a grid cell with vertices $\lambda^{(r)}\in\mathcal{G}_{\Lambda_\delta}$, $r=1,\ldots,r_g$, where the maximum distance between two points in $Q(h)$ is bounded by $\sqrt{p-1}\,h$. For $\lambda\in Q(h)$, there exist scalars $\omega_r\geq 0$ with $\sum_{r=1}^{r_g}\omega_r=1$ such that $\lambda=\sum_{r=1}^{r_g}\omega_r\lambda^{(r)}$. Given approximations $x^{(r)}$ of $x(\lambda^{(r)})$ at the vertices, we consider the convex combination $\hat{x}(\lambda)=\sum_{r=1}^{r_g}\omega_r x^{(r)}$.

\begin{lemma}
There exists a constant $K>0$ such that
$$
\|x(\lambda)-\hat{x}(\lambda)\| \leq K\sqrt{p-1}\,h+ \max_{r=1,\ldots,r_g} \|x(\lambda^{(r)})-x^{(r)}\|.
$$
\end{lemma}

\begin{proof}
As $x(\lambda)$ depends continuously on $ \lambda$ in a fixed, compact region $Q$ containing $Q(h)$, the image of this region under $x( \cdot)$ is compact. Therefore, the right hand side of \eqref{eq:first_order_opt_cond2} is uniformly bounded on $Q$. Thus, there is a constant $\max_{\mu\in Q}\|D_\lambda x(\mu)\|\leq K < \infty$. By the mean value theorem applied to the convex set $Q(h)$, we then see that $ \|x(\lambda)-x(\lambda^{(r)})\|\leq K\sqrt{p-1}\,h $. Using $\sum_{r=1}^{r_g}\omega_r=1$, we obtain
\begin{align*}
\|x(\lambda)-\hat{x}(\lambda)\|
&=\left\|\sum_{r=1}^{r_g}\omega_r
\bigl(x(\lambda)-x(\lambda^{(r)})
+x(\lambda^{(r)})-x^{(r)}\bigr)\right\|\\
&\leq K\sqrt{p-1}\,h+
\max_{r=1,\ldots,r_g}
\|x(\lambda^{(r)})-x^{(r)}\|.
\end{align*}
\end{proof}

\section{Numerical Experiments} \label{sec:NumExperiments}

In this section, we numerically investigate the proposed Pareto front tracing approach. 
We first introduce the considered classes of test problems. Next, we describe the experimental setup, including system specifications and numerical implementation details, before we discuss the numerical results.

\subsection{Example Classes}

\paragraph{Convex Quadratic functions}

First, we consider convex quadratic objective functions of the form $f_i = \frac{1}{2} (x-\chi_i)^{\top}Q_i(x-\chi_i)$ for $i \in \{1, \dots , 3\}$, where $x \in \R^n$ with $n \in \{10,20,40,80\}$. For each dimension, ten problem instances are generated. For each instance, $Q_i=A_iA_i^\top$, where the entries of $A_i\in\mathbb{R}^{n\times n}$ and $\chi_i\in\mathbb{R}^n$ are sampled independently from the standard normal distribution, and $A_i$ is resampled if singular. By construction $Q_i$ is symmetric and positive definite. Hence, each objective function $f_i$, and thus $f_\lambda$, is strongly convex and therefore strictly convex. The corresponding Pareto front is convex and connected. Moreover, $\nabla_x^2 f_\lambda(x)$ is positive definite and therefore invertible for all $\lambda \in \Lambda_\delta$. Thus, the assumptions of Theorem~\ref{thm:pareto_tracing} for tracing the Pareto front over the entire $\delta$-weight space are satisfied. As Pareto optimal solutions $x(\lambda)$ can be derived explicitly for this problem class, we can compare numerically approximated solutions by the tracing strategy directly with corresponding exact solutions. Let $Q_\lambda=\sum_{i=1}^3\lambda_iQ_i$ and $b_\lambda=\sum_{i=1}^3\lambda_iQ_i\chi_i$. The first-order optimality conditions can be rewritten as
\begin{align} \label{eq:quadr_first_order_opt_condition}
     Q_\lambda \, x(\lambda) - b_\lambda = 0 \quad
    \Leftrightarrow \quad  x(\lambda) = Q_\lambda^{-1} \, b_\lambda, 
\end{align}
For quadratic objective functions, the PDE system~\eqref{eq:first_order_opt_cond2} is given by
\begin{equation}\label{eq:first_order_opt_cond_quad}
    \partial_{\lambda_i}x(\lambda) = - Q_\lambda^{-1} \left(Q_i (x(\lambda) - \chi_i)- Q_p (x(\lambda) - \chi_p)\right),\qquad i=1,2.
\end{equation}
Indeed, by Theorem~\ref{thm:pareto_tracing}, $x(\lambda)=Q_\lambda^{-1}b_\lambda$ is the unique solution of \eqref{eq:first_order_opt_cond_quad} with initial value $x(\lambda^0)$. By Theorem~\ref{thm:pareto_tracing_init_perturbed}, the unique solution corresponding to an approximated initial value also exists over the entire $\delta$-weight space and can be derived explicitly. Hence, the dependence on initial conditions can be analyzed exactly as well. Let $x^0 = x(\lambda^0)$ denote the Pareto optimal solution with respect to the initial weight vector $\lambda^0$ and let $\widetilde{x}^0$ denote an approximate solution as in Section~\ref{subsec:InitCond}. Set $ c=\nabla_x f_{\lambda^0}(\widetilde{x}^0)$. The solution $\widetilde{x}$ corresponding to $\widetilde{x}^0$ can be derived from the first-order optimality conditions of the perturbed weighted objective $f_\lambda^c$:
\begin{align}
 \nabla_x f_\lambda^c(\widetilde{x}(\lambda)) =0 \quad \Leftrightarrow \quad \widetilde{x}(\lambda)=Q_\lambda^{-1}(b_\lambda+c). \notag
\end{align}
As $\nabla_x f_{\lambda^0}(x^0)=0$, it holds that 
$
c=\nabla_x f_{\lambda^0}(\widetilde{x}^0)-\nabla_x f_{\lambda^0}(x^0)=Q_{\lambda^0}\bigl(\widetilde{x}^0-x^0\bigr).
$
Thus, 
$$ \widetilde{x}(\lambda) = Q_\lambda^{-1} \left( b_\lambda + Q_{\lambda^0} \left(\widetilde{x}^0-x^0\right) \right).
$$
Using $x(\lambda)=Q_\lambda^{-1}b_\lambda$, the  error caused by the initial perturbation is given by
\begin{align}\label{eq:quad_fct_error_prop_exact}
\left\lVert\widetilde{x}(\lambda)-x(\lambda)\right\rVert_\infty=\left\lVert Q_\lambda^{-1}Q_{\lambda^0}\bigl(\widetilde{x}^0-x^0\bigr)\right\rVert_\infty.
\end{align}

\paragraph{WFG Toolkit}

The WFG Toolkit~\citep{hub:revtoolkit:2006} provides a framework for constructing scalable multi-objective test problems with prescribed characteristics, 
which allows us to construct two tri-objective test problems with different Pareto front shapes. The first problem has a convex Pareto front, whereas the Pareto front of the second problem has both convex and non-convex parts. For both problems, the shape functions and scaling constants of the WFG toolkit are chosen such that $\mathcal{PF}$ is explicitly known. More precisely, the underlying distance parameter is fixed at zero, while the two remaining position parameters $u_i \in [0,1], i = 1,2$, parameterize $\mathcal{PF}$~\citep{hub:revtoolkit:2006}.

The two test problems $\mathrm{WFG}_c$ and
$\mathrm{WFG}_{nc}$ are defined as
$\min_{u\in[0,1]^2}(f_1(u),f_2(u),f_{3}(u))$,
where $f_3=f_{3,c}$ and $f_3=f_{3,nc}$, respectively, with
\begin{align*}
f_1(u) &= 2(1-\cos(\tfrac{\pi}{2}u_1))(1-\cos(\tfrac{\pi}{2}u_2))
& \! f_2(u) &= 4(1-\cos(\tfrac{\pi}{2}u_1))(1-\sin(\tfrac{\pi}{2}u_2))\\
f_{3,c}(u) &= 6(1-\sin(\tfrac{\pi}{2}u_1))
& \! f_{3,nc}(u) &= 6(1-u_1-\tfrac{1}{10\pi}\cos(10\pi u_1+\tfrac{\pi}{2}))
\end{align*}
As the tracing approach is formulated for unconstrained \eqref{eq:MOP}, we define $\hat{f}_i(x) \coloneqq f_i(\Psi(x))$, where $\Psi \colon \R^2\to(0,1)^2$ is given by $\Psi(x)=\left((1+\mathrm e^{-x_1})^{-1},(1+\mathrm e^{-x_2})^{-1}\right)^\top$. As $\Psi$ is a smooth bijection between $\R^2$ and $(0,1)^2$, the reparameterization preserves Pareto dominance and, at stationary points, invertibility and positive definiteness of weighted Hessians. Note that, excluded points on the boundary of $[0,1]^2$ are not relevant in our setting, as solutions corresponding to $\lambda\in\Lambda_\delta$ lie in the interior of $\mathcal{PF}$.

In the first example, the proposed tracing approach is not limited by the Pareto front geometry, as $\mathcal{PF}$ is convex. 

While the assumptions for local tracing in Theorem~\ref{thm:pareto_tracing_neighborhood} are satisfied, strict convexity of $\hat f_\lambda$ does not hold for every $\lambda\in\Lambda$, and hence Theorem~\ref{thm:pareto_tracing} does not apply. Thus, this example allows us to investigate the tracing approach beyond the sufficient conditions for global Pareto tracing. Additionally, Pareto optimal solutions corresponding to weight vectors $\lambda \in \Lambda_\delta$ can be explicitly determined. For the box-constrained formulation, the solution mapping is given by {\small \begin{align}\label{eq:convexWFGReferenceSol}
    u(\lambda)=\left(\frac{2}{\pi}\arctan\left(\frac{3(1-\lambda_1-\lambda_2)}{\lambda_1+2\lambda_2-\sqrt{\lambda_1^2+4\lambda_2^2}}\right),\frac{2}{\pi}\arctan\left(\frac{2\lambda_2}{\lambda_1}\right)\right)^\top.
\end{align}}
Since $u(\lambda) = \Psi(x(\lambda))$, the corresponding solution of the unconstrained reformulation is given by $x(\lambda) = \Psi^{-1}(u(\lambda))$. 

In contrast to the first WFG example, the Pareto front of the second test problem also contains non-convex parts, such that the weighted sum method does not guarantee completeness. Consequently, we cannot expect to determine all non-dominated points by Pareto front tracing. 
Instead, we investigate to what extent parts of the Pareto front can be traced by numerical integration when suitable initial values are chosen. In this setting, we use the local tracing results from Theorem~\ref{thm:pareto_tracing_neighborhood}. Note that, unlike in the first problem, a 
reference solution $u(\lambda)$ is not available for the second problem.

\subsection{Experimental Setup}\label{subsec:expsetup}

Algorithm~\ref{alg:tracing_tri_obj} is implemented in Python, using numerical ODE solvers for the initial-value problems arising during the tracing procedure and nonlinear optimization solvers for weighted sum scalarizations. Unless stated otherwise, the following numerical settings and evaluation criteria are used throughout the experiments. Experiment-specific parameters are specified with the respective tests. All computations are conducted on a computer equipped with an Intel\textregistered\ Core\texttrademark i7-8700 CPU at 3.20\,GHz and 32\,GB of RAM, running a Linux operating system. The implementation uses Python~3.13.13 with NumPy~2.5.2 and SciPy~1.18.0.

\subsubsection{Numerical Implementation and Settings}\label{subsubsec:NumImplementation}

The initial-value problems in Algorithm~\ref{alg:tracing_tri_obj} are solved using the adaptive RK45 method provided by \texttt{scipy.integrate.solve\_ivp}, with tolerances $\mathrm{rtol}=10^{-3}$ and $\mathrm{atol}=10^{-6}$. The initial step size is selected automatically and no additional upper bound on the internal step size is imposed. Coordinate integrations are terminated whenever the weighted Hessian becomes nearly singular. More precisely, we use the stopping condition $ \kappa_2\bigl(H_\lambda(x)\bigr)^{-1} = {\sigma_{\min}\bigl(H_\lambda(x)\bigr)} / {\sigma_{\max}\bigl(H_\lambda(x)\bigr)} \leq 10^{-7}$. 
Approximations at weight-grid points with prescribed grid spacing $h_{\mathrm{grid}}$ are obtained by interpolation using the \texttt{t\_eval} option. All remaining solver options are left at their default values. Each tracing procedure is limited to 1800 seconds,  where the time required to compute $x^0$ is excluded. Whenever weighted sum scalarizations are solved, either to generate initial solutions or for comparison, we use the L-BFGS-B method provided by \texttt{scipy.optimize.minimize}. Analytical gradients are supplied to the solver. All remaining L-BFGS-B options are left at the default values. 

Unless stated otherwise, we set $\delta=h_{\mathrm{grid}}=0.005$ and initialize the tracing procedure with $\lambda^0=(\delta,\delta,1-2\delta)$ and a corresponding exact or numerically approximated weighted sum solution $x^0$. For the quadratic test problems, whenever an approximated initial solution is required, the optimization is started from a reproducible initial vector whose components are sampled independently from the uniform distribution on $[0,1)$.

\subsubsection{Evaluation Criteria} \label{subsubsec:evalCrit}

For each approximation $\widetilde{x}(\lambda)$ computed by the tracing strategy, we evaluate the first-order optimality residual and, if exact reference solutions $x(\lambda)$ are available, also the approximation error in the decision space by
$$
e_{\mathrm{stat}}(\lambda)
=
\lVert
\nabla_x f_\lambda\bigl(\widetilde{x}(\lambda)\bigr)
\rVert_\infty
\qquad \text{and} \qquad 
e(\lambda)
=
\lVert
\widetilde{x}(\lambda)-x(\lambda)
\rVert_\infty.
$$
Smaller values of $e_{\mathrm{stat}}$ indicate a smaller violation of the first-order optimality conditions and thus provide a useful measure when no exact reference solution is available.

We distinguish between approximations $\widetilde{x}_{\mathrm{ex}}(\lambda)$ obtained from an exact initial solution, and approximations $\widetilde{x}_{\tau}(\lambda)$ obtained from an approximated initial solution with a prescribed stationarity tolerance $\tau$. The respective stationarity residuals are
$$
e_{\mathrm{stat,ex}}(\lambda) = \lVert \nabla_x f_\lambda \bigl(\widetilde{x}_{\mathrm{ex}}(\lambda)\bigr) \rVert_\infty
\qquad \text{and} \qquad 
e_{\mathrm{stat},\tau}(\lambda) = \lVert \nabla_x f_\lambda\bigl(\widetilde{x}_{\tau}(\lambda)\bigr) \rVert_\infty.
$$
The approximation errors $e_{\mathrm{ex}}(\lambda)$ and $e_{\tau}(\lambda)$ are defined analogously. 

For each problem instance, the pointwise error measures are summarized by their maximum over all grid points at which an approximation is computed. For the quadratic test problems, these instance-wise maximum values are subsequently averaged over the ten problem instances corresponding to each decision-space dimension $n$. The resulting mean values are denoted by
$e_{\mathrm{ex}}$,
$e_{\tau}$,
$e_{\mathrm{stat,ex}}$, and
$e_{\mathrm{stat},\tau}$, respectively.

\subsection{Experiments}

We evaluate the tracing strategy on the previously introduced example classes and compare it with the Pareto Tracer~\citep{mart:predcorr:2018}.

\subsubsection{Quadratic Test Problems}

We investigate 
runtime and the dependence of initial approximation quality on the tracing accuracy. Additional results are available in the corresponding repository \citep{loehken:repository:2026}.

\begin{table}[tbp]
\centering
\footnotesize
\setlength{\tabcolsep}{2pt}

\begin{minipage}[t]{0.48\textwidth}
\centering
\vspace{0pt}
\caption{Runtime comparison of Pareto tracing and weighted sum method.}
\label{tab:time_comparison}

\medskip

\begin{tabular}{@{}rlrrrr@{}}
\toprule
$n$ & Method & $t$ [s] & $e$ & $e_{\mathrm{stat}}$ & $S$ \\
\midrule
10 & Tracing & 3.16 & 0.00526 & 0.02483 & 15.35 \\
   & WS & 48.47 & 0.00733 & 0.02138 & -- \\
   & & & & & \\
\addlinespace
20 & Tracing & 3.58 & 0.00488 & 0.04301 & 16.14 \\
   & WS & 57.88 & 0.00618 & 0.03504 & -- \\
   & & & & & \\
\addlinespace
40 & Tracing & 5.21 & 0.00512 & 0.07226 & 13.08 \\
   & WS & 68.14 & 0.00621 & 0.06354 & -- \\
   & & & & & \\
\addlinespace
80 & Tracing & 11.63 & 0.00515 & 0.14826 & 7.38 \\
   & WS & 85.77 & 0.00554 & 0.12239 & -- \\
   & & & & & \\
\bottomrule
\end{tabular}
\end{minipage}
\hfill
\begin{minipage}[t]{0.48\textwidth}
\centering
\vspace{0pt}
\caption{Influence of the initial approximation quality on tracing accuracy.}
\label{tab:init_sensitivity}

\medskip

\begin{tabular}{@{}rrrrr@{}}
\toprule
$n$ & $\tau$ & $e_\tau$
& $e_{\mathrm{init},\tau}$
& $e_{\mathrm{init,ex},\tau}$ \\
\midrule
10 & $10^{-3}$ & 0.00545
& $2.71\cdot10^{-3}$ & $2.36\cdot10^{-3}$ \\
   & $10^{-5}$ & 0.00526
& $2.12\cdot10^{-5}$ & $2.12\cdot10^{-5}$ \\
   & $10^{-7}$ & 0.00526
& $2.35\cdot10^{-7}$ & $1.92\cdot10^{-7}$ \\
\addlinespace
20 & $10^{-3}$ & 0.00499
& $2.26\cdot10^{-3}$ & $1.22\cdot10^{-3}$ \\
   & $10^{-5}$ & 0.00489
& $6.43\cdot10^{-5}$ & $1.04\cdot10^{-5}$ \\
   & $10^{-7}$ & 0.00488
& $3.02\cdot10^{-7}$ & $9.65\cdot10^{-8}$ \\
\addlinespace
40 & $10^{-3}$ & 0.00515
& $1.66\cdot10^{-3}$ & $6.86\cdot10^{-4}$ \\
   & $10^{-5}$ & 0.00512
& $2.58\cdot10^{-5}$ & $7.73\cdot10^{-6}$ \\
   & $10^{-7}$ & 0.00512
& $2.72\cdot10^{-7}$ & $1.19\cdot10^{-7}$ \\
\addlinespace
80 & $10^{-3}$ & 0.00527
& $1.93\cdot10^{-3}$ & $3.58\cdot10^{-4}$ \\
   & $10^{-5}$ & 0.00515
& $2.35\cdot10^{-4}$ & $3.50\cdot10^{-6}$ \\
   & $10^{-7}$ & 0.00515
& $8.11\cdot10^{-7}$ & $1.44\cdot10^{-7}$ \\
\bottomrule
\end{tabular}
\end{minipage}

\end{table}

\paragraph{Runtime Comparison}

We first compute approximations $\widetilde{x}_{\mathrm{ex}}(\lambda)$ using the Pareto front tracing strategy from the exact solution $x(\lambda^0)$ given by~\eqref{eq:quadr_first_order_opt_condition}. For each $\lambda \in \mathcal{G}_{\Lambda_\delta}$, we record $e_{\mathrm{stat}}$, which is then used as the prescribed stationarity tolerance for~\eqref{eq:WS}. A lower bound of $10^{-9}$ is applied to the prescribed tolerances. We denote the resulting weighted sum approximations by $\widetilde{x}_{\mathrm{WS}}(\lambda)$. 

The weighted sum scalarizations are solved sequentially over $\mathcal{G}_{\Lambda_\delta}$. We first increase $\lambda_1$ for fixed $\lambda_2$ before updating $\lambda_2$. Each scalarization is initialized with a preceding weighted sum solution, where \eqref{eq:WS} for $\lambda^{0,k+1}$ is initialized with the solution obtained at $\lambda^{0,k}$. The weighted sum approach is limited to 1800 seconds per problem instance, excluding the computation time of the exact initial solution $x^0$.

Table~\ref{tab:time_comparison} shows mean runtimes in seconds, averaged over the ten instances for each dimension. The speedup achieved by Pareto tracing is defined as $S=\overline{t}_{\mathrm{WS}}/\overline{t}_{\mathrm{trace}}$. For $\widetilde{x}_{\mathrm{WS}}$, we denote the instance-averaged maximum error measures by $e_{\mathrm{WS}}$ and $e_{\mathrm{stat,WS}}$. Since $e_{\mathrm{stat,WS}}$ is slightly smaller than $e_{\mathrm{stat,ex}}$ in all cases and $e_{\mathrm{ex}}$ and $e_{\mathrm{WS}}$ remain in a similar range, the runtime comparison is made at a comparable level of approximation accuracy. Most importantly, Pareto front tracing requires substantially less computational time than repeatedly solving \eqref{eq:WS}. The resulting speedup ranges from $7.38$ for $n=80$ to $16.14$ for $n=20$. 

\paragraph{Dependence of Approximation Quality of the Initial Value}

In this experiment, we replace the exact initial solution $x^0$ by an approximation obtained using L-BFGS-B with prescribed stationarity tolerances $\tau\in\{10^{-3},10^{-5},10^{-7}\}$. All runs are initialized from the same optimization starting point.

To quantify the effect of approximated initial values, we consider both the deviation between numerical approximations obtained from exact and perturbed initial values as well as the exact error caused by initial perturbations~\eqref{eq:quad_fct_error_prop_exact} by
$$
e_{\mathrm{init},\tau}(\lambda)
=\|\widetilde{x}_{\tau}(\lambda)-\widetilde{x}_{\mathrm{ex}}(\lambda)\|_\infty
\quad \text{and} \quad 
e_{\mathrm{init,ex},\tau}(\lambda)
=\left\|Q_\lambda^{-1}Q_{\lambda^0}
(\widetilde{x}_{\tau}(\lambda^0)-x(\lambda^0))\right\|_\infty.
$$
Analogous to the error measures in Section~\ref{subsubsec:evalCrit}, $e_{\mathrm{init},\tau}$ and $e_{\mathrm{init,ex},\tau}$ denote the averages over the ten problem instances of the instance-wise maximum values over $\mathcal{G}_{\Lambda_\delta}$.

As shown in Table~\ref{tab:init_sensitivity}, $e_{\tau}$ remains close to $e_{\mathrm{ex}}$ (Table~\ref{tab:time_comparison}). Indeed, the difference decreases as the prescribed tolerance $\tau$ is tightened. Consistent with this observation, both $e_{\mathrm{init},\tau}$ and $e_{\mathrm{init,ex},\tau}$ decrease as $\tau$ is reduced.

While $e_{\mathrm{init},\tau}$ is generally larger than $e_{\mathrm{init,ex},\tau}$, this may be explained by the additional effect of error propagation during numerical integration. Overall, the similarities between $e_{\tau}$ and $e_{\mathrm{ex}}$ suggest that, for sufficiently accurate initializations, the overall approximation error is primarily determined by numerical tracing.

Note that, while $\tau=10^{-7}$ was not attained for every instance, the resulting initial approximations were generally more accurate than for $\tau=10^{-5}$. Similar observations hold for the mean errors over the weight grid \citep{loehken:repository:2026}.

\subsubsection{WFG Test Problems}

We investigate runtime and directional averaging for the WFG examples. For $\mathrm{WFG}_{nc}$, we only report $e_{\mathrm{stat}}$, as no exact reference solution is available. Runtimes are reported as medians over ten repetitions. For RK45, we set $\mathrm{rtol}=10^{-5}$ and $\mathrm{atol}=10^{-8}$ to improve integration accuracy.

\paragraph{WFG Convex Pareto Front}

\begin{figure}[tbp]
    \centering

    \begin{minipage}[t]{0.49\textwidth}
        \centering
        \vspace{0pt}
        \small
        \setlength{\tabcolsep}{4pt}

        \captionof{table}{Runtime and accuracy comparison of tracing and weighted sum scalarizations for $\mathrm{WFG}_c$.}
        \label{tab:wfg1-directional-tradeoff}

        \resizebox{\linewidth}{!}{
            \begin{tabular}{r@{\hspace{10pt}}lrrrr}
                \toprule
                $h_{\mathrm{grid}}$ & Method & $t$ [s] & $e$ & $e_{\mathrm{stat}}$ & $S$ \\
                \midrule
                0.01 & Tracing & 6.768 & 0.00740 & 0.00002188 & 1.065 \\
                     & WS & 7.210 & 0.01162 & 0.00001445 & -- \\
                \addlinespace
                0.005 & Tracing & 17.022 & 0.01361 & 0.00003175 & 1.666 \\
                      & WS & 28.276 & 0.02898 & 0.00002318 & -- \\
                \addlinespace
                0.0025 & Tracing & 42.689 & 0.04748 & 0.00002101 & 2.579 \\
                       & WS & 110.138 & 0.05677 & 0.00001927 & -- \\
                \bottomrule
            \end{tabular}
        }
    \end{minipage}
    \hfill 
    \begin{minipage}[t]{0.49\textwidth}
        \centering
        \vspace{0pt}

        \includegraphics[
            width=\linewidth,
            height=0.55\linewidth,
            keepaspectratio
        ]{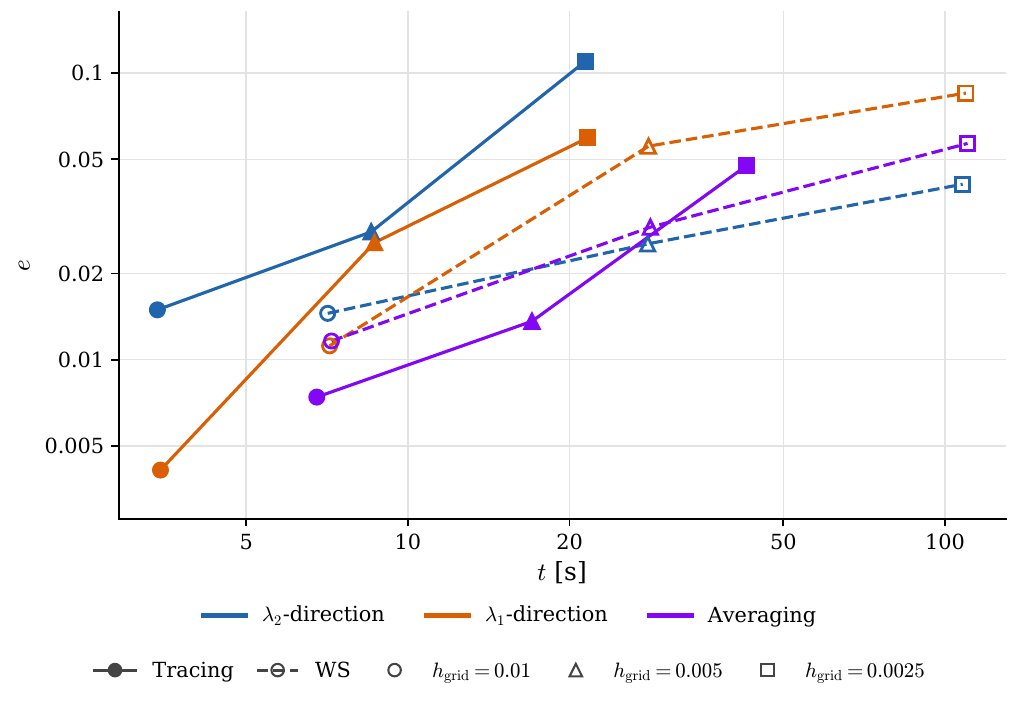}

        \captionof{figure}{Runtime and approximation error of the tracing variants for $\mathrm{WFG}_c$.}
        \label{fig:wfg_convex_comparison}
    \end{minipage}

\end{figure}

\begin{figure}[tbp]
\centering
\begin{minipage}[t]{0.49\textwidth}
    \centering
    \includegraphics[width=0.9\linewidth]{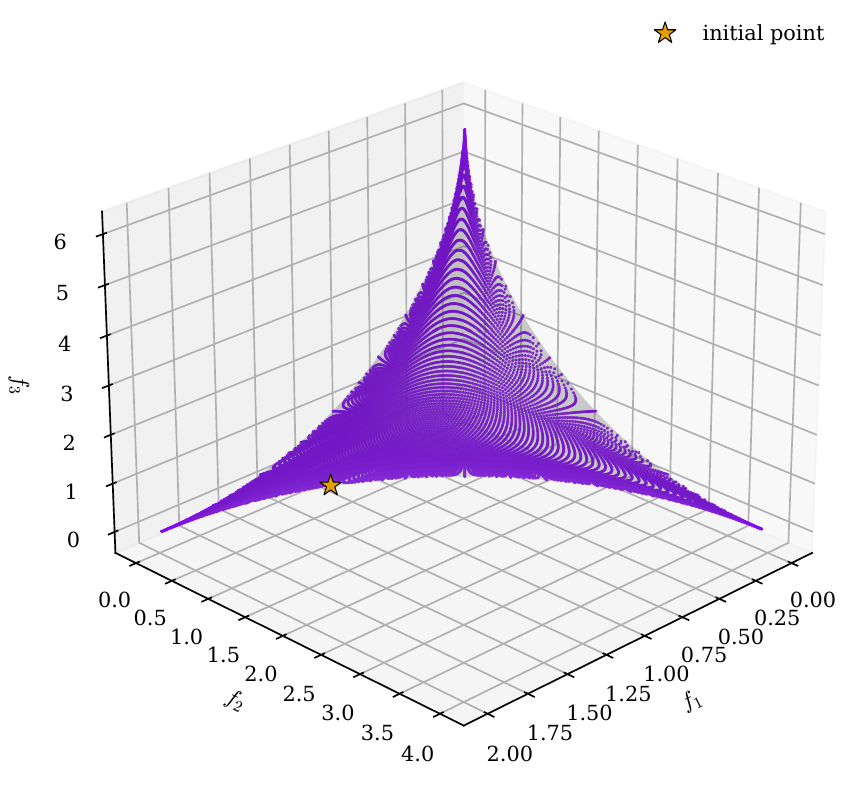}
    \caption{Pareto front coverage of $\mathrm{WFG}_{c}$.}
    \label{fig:wfg_convex_coverage}
\end{minipage}
\hfill
\begin{minipage}[t]{0.49\textwidth}
    \centering
    \includegraphics[width=0.9\linewidth]{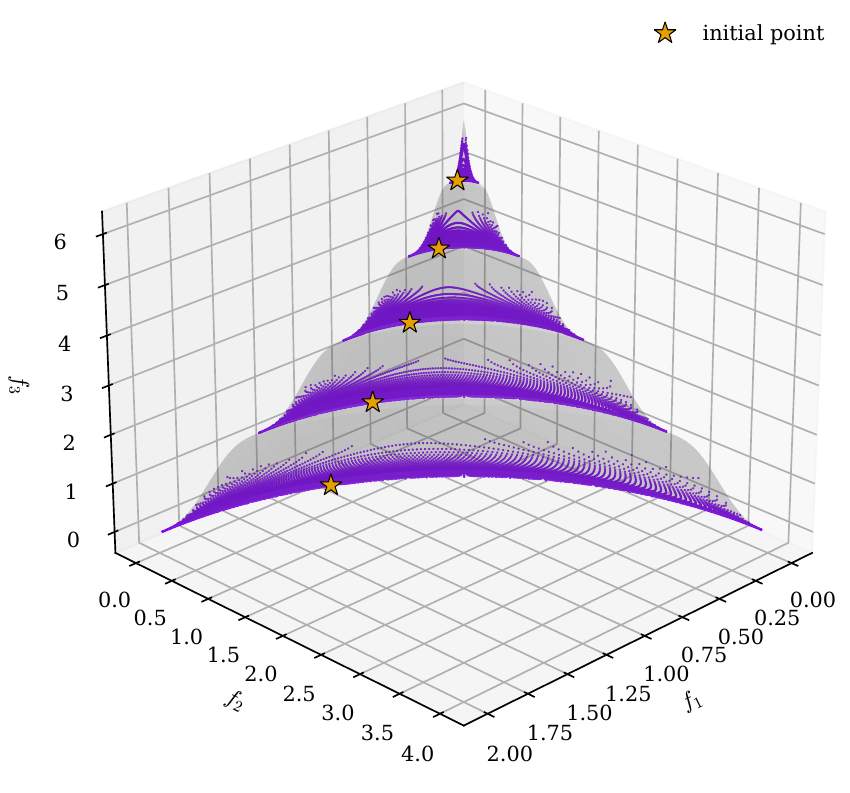}
    \caption{Pareto front coverage $\mathrm{WFG}_{nc}$.}
    \label{fig:wfg_noncvx_coverage}
\end{minipage}
\end{figure}

As illustrated in Figure~\ref{fig:wfg_convex_coverage}, the proposed tracing approach covers almost the entire $\mathcal{PF}$, except for small regions near its boundary. This indicates that its applicability extends beyond the sufficient conditions of Theorem~\ref{thm:pareto_tracing}.

Figure~\ref{fig:wfg_convex_comparison} shows that Pareto front tracing requires less runtime than repeatedly solving \eqref{eq:WS} at tested grid spacings $h_{\mathrm{grid}} \in \{0.01, 0.005, 0.0025 \}$. This holds for both directional tracing variants of Stage 2 in Algorithm~\ref{alg:tracing_tri_obj} and the averaged strategy of Stage 3. Moreover, the relative runtime advantage increases as the grid is refined (Table~\ref{tab:wfg1-directional-tradeoff}), with speedups from $S=1.065$ to $S=2.579$ for the averaged strategy. This is, however, to be expected since Pareto tracing generates approximations for most $\lambda \in \mathcal{G}_{\Lambda_\delta}$ by interpolation, whereas \eqref{eq:WS} must be solved individually for each weight vector. We note that all \eqref{eq:WS} attain their prescribed stationarity tolerances.

Directional averaging requires about twice the runtime of either directional tracing variant but yields smaller values of $e$ and $e_{\mathrm{stat}}$ in almost all cases, as also observed for the quadratic test problems \citep{loehken:repository:2026}.

In contrast to the quadratic test problems, the error $e$ seems to increase as $h_{\mathrm{grid}}$ decreases (Table~\ref{tab:wfg1-directional-tradeoff}). However, further numerical experiments show that this is caused by additional weight vectors included in the finer grids, particularly those closer to the boundary of the $\delta$-weight space. When the averaged approximations obtained for the finer grids are evaluated only at weight vectors corresponding to $h_{\mathrm{grid}}=0.01$, $e$ decreases from $0.007404$ to $0.003746$ and $0.003712$, respectively.

\paragraph{WFG Non-Convex Pareto Front} As we cannot expect to trace the entire $\mathcal{PF}$ from a single initial solution, we apply Algorithm~\ref{alg:tracing_tri_obj} with different initial values, each corresponding to a different locally convex region of $\mathcal{PF}$. Those are obtained by solving \eqref{eq:WS} with manually varied box constraints. Figure~\ref{fig:wfg_noncvx_coverage} shows that the proposed approach successfully traces large regions of $\mathcal{PF}$, including local regions that cannot be covered by unconstrained weighted sum scalarizations.

For the runtime comparison, we solve \eqref{eq:WS} with additional box constraints to restrict the optimization to the respective local regions. Table~\ref{tab:wfg2_results} shows that, for $h_{\mathrm{grid}}=0.01$, the averaged tracing variant requires more computational time than repeatedly solving corresponding~\eqref{eq:WS}. However, the median speedup increases to $S=1.091$ and $S=1.905$ as the grid is refined. Moreover, as observed for the convex WFG example, following only one integration direction reduces the runtime by approximately half, providing a faster alternative when computational time is critical. Importantly, however, even when the weighted sum approach is faster, numerical tracing only requires one suitable initial solution for each local region, whereas the weighted sum approach additionally requires appropriately chosen box constraints, whose determination may be computationally challenging.

The influence of the integration direction in Stage~2 of Algorithm~\ref{alg:tracing_tri_obj}, as well as the subsequent averaging in Stage~3, becomes particularly interesting, as we can only expect to trace the Pareto front locally. For each variant, the percentage of reached grid points ($r_{\mathcal G}$) and the maximum stationarity residual $e_{\mathrm{stat}}$ are evaluated over all five tracing procedures. For the error evaluation, we further distinguish between all points reached by each variant ($e_{\mathrm{stat,all}}$) and the common points reached by all three variants ($e_{\mathrm{stat,com}}$). Interestingly, integration in $\lambda_1$-direction yields smaller maximum stationarity residuals than both $\lambda_2$-direction and averaging when evaluated for common grid points. Thus, in contrast to the previous test problems, averaging does not necessarily improve accuracy. However, combining both integration directions increases the overall coverage, as the directional variants do not always reach the same grid points.

\begin{table}[tbp]
\centering
\caption{Runtime, percentage of reached grid points, and stationarity residuals for $\mathrm{WFG}_{nc}$ over five
initializations. Runtimes and speedups are medians over ten repetitions.}
\label{tab:wfg2_results}

\footnotesize
\setlength{\tabcolsep}{3pt}

\resizebox{\textwidth}{!}{
\begin{tabular}{@{}llrrrrrrrr@{}}
\toprule
$h_{\mathrm{grid}}$ & Variant
& $t$ [s] & $t_{\mathrm{WS}}$ [s] & $S$
& $r_{\mathcal G}$ [\%]
& $e_{\mathrm{stat,all}}$
& $e_{\mathrm{stat,com}}$
& $e_{\mathrm{stat,WS,all}}$
& $e_{\mathrm{stat,WS,com}}$ \\
\midrule
0.01
& $\lambda_1$-direction
& 36.84 & 45.48 & 1.234
& 92.67 & 0.0958 & 0.0000791 & 0.0844 & 0.0366 \\
& $\lambda_2$-direction
& 36.56 & 44.74 & 1.223
& 92.78 & 0.000365 & 0.000365 & 0.0366 & 0.0366 \\
& Averaging
& 72.94 & 45.60 & 0.625
& 93.00 & 0.0958 & 0.000180 & 0.0844 & 0.0366 \\
\midrule
0.005
& $\lambda_1$-direction
& 82.39 & 178.05 & 2.161
& 92.34 & 0.000228 & 0.000228 & 0.0141 & 0.0141 \\
& $\lambda_2$-direction
& 83.10 & 176.61 & 2.123
& 92.41 & 0.0141 & 0.0141 & 0.0141 & 0.0141 \\
& Averaging
& 164.59 & 179.45 & 1.091
& 92.67 & 0.00676 & 0.00676 & 0.0141 & 0.0141 \\
\midrule
0.0025
& $\lambda_1$-direction
& 181.83 & 684.74 & 3.761
& 91.33 & 0.000146 & 0.000146 & 0.00946 & 0.00946 \\
& $\lambda_2$-direction
& 182.03 & 673.04 & 3.689
& 92.07 & 0.0809 & 0.00399 & 0.0787 & 0.000718 \\
& Averaging
& 361.03 & 688.20 & 1.905
& 92.45 & 0.0809 & 0.00224 & 0.0787 & 0.000718 \\
\bottomrule
\end{tabular}
}
\end{table}

\subsubsection{Comparison with Predictor Corrector Methods}

To compare the proposed tracing strategy with the Pareto Tracer (PT) introduced by \cite{mart:predcorr:2018}, we investigate computational cost and approximation accuracy on $\mathrm{WFG}_{c}$ at comparable levels of representation quality of generated outcome sets, measured in terms of coverage~\citep{sayin:repr:2000}. Let $A=\{f(x):x\in \mathcal{X}_{\mathrm{approx}}\} \subset \R^p$ denote a generated outcome set and $\mathcal{Y}_R \subset \R^p$ a reference set accurately describing $\mathcal{PF}$. Using the weighted Tchebycheff distance, the coverage error can be computed as
$$
\operatorname{CE}(A, \mathcal{Y}_R)=\max _{y \in \mathcal{Y}_R} \min_{a \in A} \max_{i=1, \ldots, p} \left\{w_i \bigl|y_i-a_i\bigr|\right\} ,
$$
where $w_i={1}/{(y_i^{\max}-y_i^{\min})}$, $i=1,\ldots,p$, and $y_i^{\max}$ and $y_i^{\min}$ are the maximum and minimum of $f_i$ on $\mathcal{PF}$, respectively. For $\mathrm{WFG}_{c}$, these ranges are explicitly known. 

Using the analytical solution in~\eqref{eq:convexWFGReferenceSol}, we construct a reference set $\mathcal{Y}_R$ by uniform sampling on $[0,1]^2$, yielding $9,878,449$ points. Due to bounded partial derivatives of the normalized objectives, $\operatorname{CE}(\mathcal{Y}_R, \mathcal{PF})$ with respect to $\mathcal{PF}$ is bounded by $10^{-3}$. To avoid increasing the runtime of PT, particularly by covering boundary regions of $\mathcal{PF}$ that may be of limited interest, we restrict the comparison to an interior region by selecting only reference points whose corresponding weight vectors satisfy $\lambda_i\geq 0.1$, $i=1,2,3$, yielding $3,450,567$ reference points.

For both approaches, we select configurations that yield comparable coverage errors. For tracing, we set $\delta=0.01$ and vary $h_{\mathrm{grid}}$, while for PT, we vary $\mathrm{PCStepObj}$, which defines the desired Euclidean distance between consecutive objective vectors. All remaining settings of PT are left at their defaults. Coverage errors are computed over all generated objective vectors, while $e$ and $e_{\mathrm{stat}}$ are evaluated only for points corresponding to the considered interior part of $\mathcal{Y}_R$, where $N_{\mathrm{acc}}$ denotes the number of evaluated points. Table~\ref{tab:wfg1_ParetoTracer_comparison} presents median runtimes over ten test runs.

Indeed, the Pareto Tracer yields approximations of Pareto optimal solutions 
with $e$ and $e_{\mathrm{stat}}$ considerably smaller. However, the proposed tracing strategy attains comparable interior coverage with substantially smaller output sets, which is a large advantage for a decision maker, who usually prefers smaller sets of solutions to choose from. Moreover, the tracing strategy requires only about $3.5$--$3.7\%$ of the runtime needed by PT. 

\begin{table}[tbh]
    \centering
    \small
    \setlength{\tabcolsep}{4pt}
    \begin{tabular}{llrrrrrr}
        \toprule
        Method & Parameters & $|A|$ & $e_{\mathrm{cov}}$
        & $t$ [s] & $N_{\mathrm{acc}}$ & $e$ & $e_{\mathrm{stat}}$ \\
        \midrule
        Tracing
        & $h_\mathrm{grid} = 0.97/6$
        & 28 & 0.20068 & 0.219 & 10
        & $4.43 \cdot 10^{-5}$ & $4.93 \cdot 10^{-6}$ \\

        PT
        & $\mathrm{PCStepObj} = 0.15$
        & 435 & 0.19667 & 5.964 & 212
        & $7.11 \cdot 10^{-15}$ & $5.00 \cdot 10^{-16}$ \\

        \midrule
        Tracing
        & $h_\mathrm{grid} = 0.97/8$
        & 45 & 0.11251 & 0.293 & 21
        & $8.82 \cdot 10^{-5}$ & $9.29 \cdot 10^{-6}$ \\

        PT
        & $\mathrm{PCStepObj} = 0.10$
        & 842 & 0.12531 & 8.326 & 433
        & $9.33 \cdot 10^{-15}$ & $4.16 \cdot 10^{-16}$ \\
        \bottomrule
    \end{tabular}
    \caption{Runtime and accuracy comparison between the proposed tracing strategy and the Pareto Tracer at comparable coverage values on $\mathrm{WFG}_{c}$.}
    \label{tab:wfg1_ParetoTracer_comparison}
\end{table}

\section{Conclusions}\label{sec:Conclusion}

In this work, we propose a novel approach for tracing Pareto optimal solutions of unconstrained multi-objective optimization problems with an arbitrary number of objectives by numerical integration. Based on a PDE system describing Pareto critical solutions, that is derived from first-order optimality conditions of weighted sum scalarizations, we state conditions for local existence, uniqueness, and Pareto optimality of its solutions. Under additional assumptions, we show that the solution can be extended to the entire weight space, such that we can trace the entire interior of the Pareto front. Finally, we develop a tracing strategy for numerically solving the PDE system and provide a comprehensive error analysis.

In our numerical experiments, we observe significant time savings compared with repeatedly solving weighted sum scalarizations. Compared with predictor–corrector methods, the proposed strategy yields similar representation quality with respect to coverage, requiring considerably fewer points and lower computational costs. 
Notably, while the weighted sum method guarantees completeness only for convex \eqref{eq:MOP}, our approach is also applicable to Pareto fronts with non-convex parts.

Particularly for problems with Pareto fronts containing both convex and non-convex parts, the proposed tracing strategy relies on suitable initial solutions. Therefore, future research includes developing strategies for finding initial solutions in different local regions. To improve approximation accuracy, the proposed tracing strategy could be combined with predictor–corrector methods by incorporating corrector steps. Another direction includes alternative integration strategies, such as directly integrating along line segments from the initial weight vector to different weight vectors, which could be particularly interesting when only a small number of Pareto optimal solutions is desired. Further research could also address the extension to constrained problems. Moreover, this work serves as a starting point for the development of momentum methods based on port Hamiltonian control strategies to recover the Pareto fronts of nonconvex shape optimization problems.

\appendix

\section*{Acknowledgments}

Matthias Bolten, Hanno Gottschalk, and Kathrin Klamroth acknowledge funding by the Deutsche Forschungsgemeinschaft (DFG, German Research Foundation) – Project-ID 531152215 – CRC 1701. We also thank Claudia Totzeck for constructive discussions.

\noindent\textit{Use of AI tools.}
During the preparation of this manuscript, ChatGPT and DeepL were used to improve language and overall readability. Additionally, ChatGPT was used for code restructuring and optimization (e.\,g., improved error handling and performance improvements) and for exploring and discussing different approaches to mathematical proofs, both in the early stages of the research and during proofreading of mathematical arguments. The manuscript has been fully revised and verified by the authors. The authors take full responsibility for all content.

\bibliographystyle{abbrvnat}
\bibliography{literatur}

\end{document}